\pdfoutput=1
\documentclass[hidelinks,onefignum,onetabnum]{siamart251216}

\usepackage{graphicx} %

\usepackage{amsfonts}
\usepackage{stmaryrd}
\usepackage{adjustbox}
\usepackage{float}
\usepackage{booktabs}
\usepackage{multirow}
\usepackage{hyperref}
\usepackage{bm}
\usepackage{tikz}
\usepackage{tcolorbox} %
\usepackage{xcolor}
\usepackage{url}

\ifpdf
  \DeclareGraphicsExtensions{.eps,.pdf,.png,.jpg}
\else
  \DeclareGraphicsExtensions{.eps}
\fi

\DeclareMathOperator{\dive}{div}
\DeclareMathOperator{\grad}{\nabla}

\newcommand{\R}{\mathbb R}
\newcommand{\Hone}[1]{H^1(\Omega_{#1}) }

\newcommand{\sumGammaij}{\sum_{\substack{j=0\\ j \neq i}}^N}

\newcommand{\Gammaij}{\Gamma_{ij}}

\headers{Multi-GPU Acceleration of Interface Problems in FEniCS}{E. Centofanti, B. Pachev, J. Steyer, J.D. Trotter, P. Benedusi}

\title{Multi-GPU Acceleration of Interface Problems in FEniCS, with Application to Cardiac Modeling}

\author{Edoardo Centofanti\thanks{Università di Pavia, Italy} \and Benjamin Pachev\thanks{University of Massachusetts Amherst, USA} \and
Joshua Steyer\thanks{Karlsruhe Institute of Technology, Germany} \and
James D. Trotter\thanks{Simula Research Laboratory, Norway} \and
Pietro Benedusi\thanks{Università della Svizzera Italiana, Switzerland (\email{benedp@usi.ch})}}

\begin{document}

\maketitle

\begin{abstract}
GPU acceleration of finite element assembly is well established for standard formulations, but mixed-dimensional problems, in which volumetric subdomains are coupled via lower-dimensional interfaces, remain unsupported in automated frameworks such as FEniCS.
We extend cuDOLFINx, a GPU-enabled extension of the DOLFINx finite element framework, to assemble such interface problems on GPUs while preserving the high-level FEniCS programming model. Combined with GPU-resident Krylov solvers and algebraic multigrid preconditioning, this enables end-to-end GPU execution of the EMI (Extracellular–Membrane–Intracellular) model, an elliptic interface problem from cardiac modeling in which many intracellular subdomains are coupled trough membranes. We compare the GPU implementation against a CPU-based one through strong and weak scaling experiments on regular and realistic 3D cardiac geometries with up to 8\,000 cells and 22 million degrees of freedom. In terms of assembly, a single GH200 node is roughly equivalent to two 144-core CPU nodes in strong scaling, and to two to three nodes in weak scaling, depending on geometry. For the linear solver, performance depends on geometry, preconditioning, and problem scale. On the regular geometry, GPU advantages emerge mainly at larger problem sizes and node counts. On the realistic geometry, with the same algebraic multigrid preconditioner for both architectures, the GPU is faster at all tested scales. The implementation is available as open-source software and applies to other problems with interface coupling.
\end{abstract}



\section{Introduction}
Modern GPU architectures now dominate many HPC systems, offering substantially higher throughput and memory bandwidth than multicore CPUs, which makes them attractive for large-scale finite element simulations, where both assembly and linear solves are major costs. While GPU-enabled linear algebra is available in several libraries~\cite{Hypredocs,Petscdocs,heinlein2018frosch}, GPU acceleration of finite element assembly remains considerably more challenging, due to complex mesh data structures, code generation, and memory management.

FEniCS~\cite{BarattaEtal2023} has become one of the most widely adopted open-source finite element frameworks thanks to its expressive Python interface and automated code generation, enabling the rapid development of finite element applications from high-level mathematical formulations. Recently, the experimental CUDA extension cuDOLFINx~\cite{pachev2026cudolfinx} demonstrated that matrix and vector assembly can be efficiently offloaded to GPUs while requiring only minimal modifications to existing FEniCS applications. By preserving the high-level programming model of FEniCS, cuDOLFINx makes GPU acceleration accessible without requiring users to rewrite their codes using GPU-specific programming models.

Despite these advances, an important class of finite element applications remains beyond the capabilities of current GPU assembly frameworks: mixed-dimensional problems. Many multiphysics models involve interfaces, where multiple subdomains are coupled through lower-dimensional manifolds. Such problems require the simultaneous assembly of volumetric and interface contributions, posing challenges that are fundamentally different from those encountered in standard finite element formulations. Extending GPU assembly to these interface problems is therefore a necessary step towards broadening the applicability of GPU-accelerated finite element software.

A relevant example of a mixed-dimensional problem is the EMI (Extracellular–Membrane–Intracellular) or cell-by-cell model~\cite{tveito2021modeling,jaeger2021efficient}, used in computational cardiology to describe electrical propagation at cellular resolution~\cite{SteyerCheginiPotseetal.2023}. Unlike homogenized monodomain and bidomain models~\cite{franzone2014mathematical}, EMI represents intracellular and extracellular domains explicitly, coupled through membrane interfaces. Advances in micrometer-scale 3D reconstruction now make realistic cellular geometries available~\cite{shapson2024petavoxel}, but even a few thousand cells lead to $\mathcal{O}(10^7)$ degrees of freedom, making both assembly and linear solution major bottlenecks. Scalable solvers have been developed, including algebraic multigrid~\cite{centofanti2025parallel} and domain decomposition methods such as BDDC~\cite{gobel2026bddc} and GDSW~\cite{centofantigdsw}, but mostly on regular geometries or structured meshes. Realistic unstructured geometries have so far been considered mainly in computational neuroscience~\cite{benedusi2024modeling}, where only the extracellular space and a single intracellular compartment are coupled. The cardiac setting is considerably more challenging: neighboring cells are electrically connected through gap junctions, so a large number of intracellular subdomains must be coupled simultaneously, giving rise to highly complex interface problems.

To address this challenge, in this work we extend cuDOLFINx to support GPU assembly of interface problems, enabling for the first time GPU-accelerated assembly of both matrices and right-hand sides for the cardiac EMI model on large realistic unstructured meshes. Rather than introducing a problem-specific implementation, our approach generalizes the cuDOLFINx framework to interface formulations involving coupled volumetric and codimension-one domains while preserving the high-level FEniCS programming model. Combined with GPU-resident iterative solvers and preconditioners, this enables efficient end-to-end GPU execution of large-scale EMI simulations. We assess the proposed implementation through weak and strong scalability studies on state-of-the-art HPC architectures using both idealized regular geometries and realistic cardiac meshes.
We stress that this study goes beyond the cardiac application, since many other problems share a similar structure, where elliptic differential operators are coupled via jumps and flux conditions on interfaces, e.g., in fracture modeling, where these conditions are generalized via Robin–Wentzell interface conditions \cite{favino2025fractures}.

This paper is organized as follows: in Section~\ref{sec:emi} we introduce the cardiac EMI model. In Section~\ref{sec:impl} we describe the implementation, in terms of GPU acceleration and solution strategy. In Section~\ref{sec:models}, we present the computational geometries and meshes used for benchmarking. Finally, in Section~\ref{sec:results}, we report the weak and strong scalability results obtained using state-of-the-art HPC systems.

\section{The EMI model}\label{sec:emi}
In this section, we describe the cardiac EMI model, also known as cell-by-cell model \cite{benedusi2024modeling, tveito2021modeling}, using the same formulation as in \cite{huynh2023convergence}.  
We consider the domain $\Omega\in\R^3$, as composed of $N+1$ subdomains, corresponding to $N$ cells $\Omega_i$, with $i=1,\dots,N$, and an extracellular medium $\Omega_0$.
Ionic exchanges through cell membranes are responsible for the interactions between cells and the extracellular space.
Adjacent cells can be connected through gap junctions, which contain protein channels responsible for the flow and exchange of ionic species. 
From a mathematical perspective, the EMI model is intrinsically mixed-dimensional: it couples volumetric PDEs in the intra- and extracellular subdomains with lower-dimensional interface equations (typically ODEs) posed at the membranes.

The EMI model describes the time evolution of the electrical potential $u_i:\overline{\Omega}_i\times(0,T]\to\mathbb{R}$ within the subdomain $\Omega_i$ and across the interfaces $\Gamma_{ij}=\overline{\Omega}_i \cap \overline{\Omega}_j$, for $i,j=0,\ldots,N$, and has the following form:
\begin{equation}\label{eq: cell-by-cell model} 
	\begin{cases}
		- \dive (\sigma_i \grad u_i) = 0		& \text{in}\quad\Omega_i,\\
		- \bm{n}_i \cdot\sigma_i \grad u_i =  \bm{n}_j \cdot\sigma_j \grad u_j = C_\mathrm{M} \dfrac{\partial v_{ij} }{\partial t} + I_{ij}(v_{ij}, \bm{w_{ij}})		& \text{in } \Gamma_{ij}, \quad i\neq j,		\\		
		\dfrac{\partial \bm{w}_{ij}}{\partial t} - R(v_{ij}, \bm{w}_{ij}) = 0  & \text{in } \Gamma_{ij},\quad i\neq j.
	\end{cases}
\end{equation}
Here, $v_{ij} = u_i - u_j$ is the potential jump between cell $i$ and cell $j$, across the common boundary $\Gamma_{ij}$, $\sigma_i>0$ is the conductivity coefficient corresponding to the $i$th subdomain, $\bm{n}_i$ is the outward normal to $\Omega_i$, and $C_\mathrm{M}>0$ is the specific membrane capacitance per unit area of the membrane surface. We notice that here we consider a common capacitance; a more general model can include a set of $C_{ij}$, i.e. with the option of varying capacitance depending on the specific interface.
The $I_{ij}$ term (possibly non-linear in $v_{ij}$) represents the current flowing through the cellular membrane because of ionic exchanges.
Typically 
 \begin{equation}\label{eq:F_term}
	I_{ij}(v,\bm{w}_{ij}) = \left\{ 
        \begin{array}{ll}
            R_\mathrm{g}^{-1}v          & \text{if}\quad i\neq0 \quad\text{and}\quad j\neq0, \\ 
	        I_\text{ion}(v, \bm{w}_{ij}) & \text{otherwise},	        
        \end{array}
        \right.
\end{equation}
where $R_\mathrm{g}$ is the gap junction resistance, and $I_\text{ion}$ is a specific ionic model, which depends on a set of $M$ auxiliary variables $\bm{w}_{ij}:\Gamma_{ij}\to\mathbb{R}^M$ (e.g., gating variables), related to the behavior of the ionic channels in the cellular membranes. In this setting, $R$ is usually a non-linear function describing membrane dynamics and recovery. 
As a result, the second equation in \eqref{eq: cell-by-cell model} represents the physiological conservation of current fluxes between the intra- and the extracellular spaces, as well as an ODE system modeling ion flow dynamics. 

Finally, we close the EMI system \eqref{eq: cell-by-cell model} with  homogeneous Neumann boundary conditions: for $i=0,\ldots,N$
$$\bm{n}_i\cdot \grad u_i = 0 \quad \text{in}\quad\partial\Omega\cap\Omega_i,$$
and with an initial condition on membrane voltage $v_{ij}=v_{ij}^\text{in}$ for $t=0$, given the initial voltage $v_{ij}^\text{in}:\Gamma_{ij}\to\mathbb{R}$ and initial auxiliary variables $\bm{w}^\text{in}_{ij}$. According to  \eqref{eq: cell-by-cell model}, the solving potential can be determined up to a constant, so that uniqueness has to be enforced, for example with an additional point-wise Dirichlet condition $u_0(\bm{x}_\mathrm{D},t)=0$ for $\bm{x}_\mathrm{D}\in\partial\Omega_0$, or at the solution strategy level, imposing a constant\footnote{Experiments show that these two approaches are essentially equivalent}. 

We discretize system \eqref{eq: cell-by-cell model} in space using the finite element method, while we use an implicit-explicit (IMEX) strategy to discretize in time: first, the auxiliary variables $\bm{w}_{ij}$ are computed using the jump $v_{ij}$ from the previous time step with an explicit treatment of $R$; then, the ionic current $I_{ij}$ is computed. Finally, we solve the EMI equations with respect to the electric potential. 
In the following sections, for completeness, we report the finite element discretization, following the notation in \cite{benedusi2024dense}.

\subsection{Weak formulation}
    By integrating the first equation in (\ref{eq: cell-by-cell model}) by parts over each $\Omega_i$ and incorporating the second equation (defined on the membrane) into the resulting interface terms, we obtain the weak formulation for the $i$-th subdomain. Specifically, for $t \in (0, T]$, we seek $u_i \in H^1(\Omega_i)$ for $i = 0, \ldots, N$, such that:
    \begin{equation*}
        \int_{\Omega_i} \sigma_i \grad u_i \cdot \grad \phi_i \, \mathrm{d}x + \sumGammaij \int_{\Gammaij} \left( C_\mathrm{M} \dfrac{\partial v_{ij}}{\partial t} + I_{ij} \right) \phi_i \, \mathrm{d}s = 0,\qquad \forall \phi_i \in \Hone{i}.
    \end{equation*}
Substituting $v_{ij}=u_i-u_j$, we can express the previous equation as    
\begin{equation}\label{eq:weak}
        \int_{\Omega_i} \sigma_i \grad u_i \cdot \grad \phi_i \, \mathrm{d}x + \int_{\Gamma_i} C_\mathrm{M}\dfrac{\partial u_i}{\partial t} \phi_i \, \mathrm{d}s - \sumGammaij \int_{\Gammaij} C_\mathrm{M} \dfrac{\partial u_j}{\partial t} \phi_i \, \mathrm{d}s = -\sumGammaij \int_{\Gammaij}  I_{ij}  \phi_i \, \mathrm{d}s,
    \end{equation}
with $\Gamma_i=\partial\Omega_i\backslash \partial\Omega$, for all $\phi_i \in \Hone{i}$. 

\subsection{Time discretization}
We uniformly discretize the time domain $[0,T]$ with $N_t$ intervals, so that $\Delta t=T/N_t$. We describe the time integration of a single time step; the complete time evolution then follows accordingly. Given the initial conditions $v_{ij}^\text{in}$ and $\bm{w}_{ij}^\text{in}$, we discretize the last equation in \eqref{eq: cell-by-cell model} with explicit Euler, i.e., we obtaining $\bm{w}_{ij}$ in $\Gammaij$.
Then, we define $I_{ij}=I_{ij}(v_{ij}^\text{in},\bm{w}_{ij})$ and equation \eqref{eq:weak} becomes 
\begin{equation}\label{eq:weak_2}
        \tau\int_{\Omega_i} \sigma_i \grad u_i \cdot \grad \phi_i \, \mathrm{d}x +\int_{\Gamma_i}u_i\phi_i \, \mathrm{d}s -\sumGammaij \int_{\Gammaij} u_j \phi_i \, \mathrm{d}s = \sumGammaij \int_{\Gammaij}  (v_{ij}^\text{in} - \tau\,I_{ij})  \phi_i \, \mathrm{d}s,
\end{equation}
where $\tau=\Delta t /C_\mathrm{M}$ and the left hand side potentials $u_i$ and $u_j$ are evaluated at $t=\Delta t$. For $i=0,\ldots,N$, we define the bilinear and linear forms
\begin{equation} \label{eq:d_i}
        d_i (u_i, \phi_i) := \tau a_i(u_i, \phi_i) + p_i(u_i, \phi_i),\quad f_i (\phi_i) := \sumGammaij \int_{\Gammaij}  (v_{ij}^\text{in} - \tau\,I_{ij})  \phi_i \, \mathrm{d}s.
\end{equation}
where
\begin{equation*} \label{eq: a_i p_i}		
			a_i (u_i, \phi_i) := \int_{\Omega_i} \sigma_i \grad u_i \cdot \grad \phi_i \, \mathrm{d}x,					\quad 
		p_i (u_i, \phi_i) := \int_{\Gamma_i}u_i\phi_i \, \mathrm{d}s -\sumGammaij \int_{\Gammaij} u_j \phi_i \, \mathrm{d}s.		
	\end{equation*}
Given $\mathcal{H}^1=\Hone{0}\times\Hone{1}\times\ldots\times\Hone{N}$, the EMI single-time-step weak formulation reads: find $(u_0,u_1,\ldots,u_N) \in \mathcal{H}^1$  so that 
$$ \sum_{i=0}^N d_i (u_i, \phi_i) = \sum_{i=0}^N f_i (\phi_i),\quad\forall (\phi_0,\phi_1,\ldots,\phi_N) \in \mathcal{H}^1.$$

\subsection{Space discretization}\label{sec:space_disc}
For the spatial discretization, we use the finite element method with linear elements (P1) on tetrahedral meshes. For $i=0,\ldots,N$, let $V(\overline{\Omega}_i)$ be the finite element space of piecewise continuous linear functions $u_{i,h}$ in $\Omega_i$, encoded by $n_i$ degrees of freedom.
The continuous problem becomes: find $(u_{0,h},u_{1,h},\ldots,u_{N,h})\in V(\overline{\Omega}_0)\times V(\overline{\Omega}_1) \times\ldots\times V(\overline{\Omega}_N)$ such that
\begin{equation}\label{eq:problem}
    \sum_{i=0}^N d_i (u_{i,h}, \phi_{i,h}) = \sum_{i=0}^N f_i (\phi_{i,h})
\end{equation}
is satisfied for all $(\phi_{0,h},\phi_{1,h},\ldots,\phi_{N,h})\in V(\overline{\Omega}_0)\times V(\overline{\Omega}_1) \times\ldots\times V(\overline{\Omega}_N)$.
Alternatively, in matrix form, Eq.~\eqref{eq:problem} may be written as:
\begin{equation}\label{eq:monolithic}
    \mathcal{K}\bm{u} = \bm{f},\quad \text{ with } \mathcal{K} = \sum_{i=0}^N\mathcal{K}_i,\quad \mathcal{K}_i = \tau A_i + M_i,
\end{equation}
where $n=\sum_i^N n_i$ is the global number of degrees of freedom (DoFs), $A_i,M_i\in\mathbb{R}^{n\times n}$ are stiffness and mass matrices corresponding to $a_i$ and $p_i$ in Eq. \eqref{eq:d_i}, respectively, while $\bm{u},\bm{f}\in\mathbb{R}^n$ collect the finite element solution coefficients and the right and side in Eq.~\eqref{eq:problem}. 
In this setting, the mass matrices contain only interface contributions. We remark that $\mathcal{K}\in\mathbb{R}^{n\times n}$ is symmetric by construction. We refer to \cite{benedusi2024dense} %
for a more detailed description of the block structure of the matrix $\mathcal{K}$, with an analysis of the interplay between $n_i, n$ and the number of degrees of freedom of the membrane. Under uniform refinement, $\mathcal{K}$ is shown to be spectrally equivalent to $N+1$ independent Poisson problems, one for each cell.

Finally, we remark that a discontinuous Galerkin (DG) discretization over the entire domain can also resolve interface jumps, as in \cite{Ellingsrud2025}. 
This avoids restricted assembly in $\Omega_i$, but increases the number of degrees of freedom. Restricted assembly instead captures interface discontinuities while using a continuous Galerkin discretization within each cell, reducing the degrees of freedom compared to a global DG scheme.

\section{Implementation}\label{sec:impl}
For the finite element discretization, we use the FEniCS framework (DOLFINx v0.9.0)~\cite{BarattaEtal2023}, which
translates weak forms expressed in the Unified Form Language (UFL) into low-level C code via the FEniCS Form Compiler (FFCx).
PETSc~\cite{Petscdocs}, providing a parallel implementation of linear solvers. We use the multiphenicsx~\cite{multiphenicsx} library to assemble mixed-dimensional variational problems coupled across subdomains, as in~\cite{fevola2021automated, fumagalli2024discontinuous, causemann2026pulsating}. The code is openly available in the cardioEMI\footnote{\url{https://github.com/pietrobe/cardioEMI}} repository.
We first describe how mixed-dimensional weak forms are set up and assembled on CPUs, we then motivate the need for GPU-accelerated assembly and describe how it is realized in cuDOLFINx. Finally, we detail the solution strategy adopted for the resulting linear using the CUDA backend for PETSc.

\subsection{Assembly of the weak forms}\label{sec:weak_forms}
\begin{figure}
    \centering
    \begin{adjustbox}{max width=\linewidth}
    \begin{tikzpicture}[scale=3]
\def\n{10} %

\draw[thick] (0,0) rectangle (1,1);

\foreach \i in {0,...,\n} {
    \foreach \j in {0,...,\n} {
        \pgfmathsetmacro\x{\i/\n}
        \pgfmathsetmacro\y{\j/\n}
        \ifnum\i<\n
            \ifnum\j<\n
                \draw[gray!60, thin] (\x,\y) -- ({(\i+1)/\n},\y);
                \draw[gray!60, thin] (\x,\y) -- (\x,{(\j+1)/\n});
                \draw[gray!60, thin] (\x,\y) -- ({(\i+1)/\n},{(\j+1)/\n});
            \fi
        \fi
    }
}

\coordinate (A) at (0.1,0.1);
\coordinate (B) at (0.1,0.5);
\coordinate (C) at (0.5,0.5);
\coordinate (D) at (0.5,0.1);

\fill[red!30, opacity=0.5] (A) -- (B) -- (C) -- (D) -- cycle;

\draw[thick] (A) -- (B) -- (C) -- (D) -- cycle;

\coordinate (A1) at (0.2,0.5);
\coordinate (B1) at (0.2,0.8);
\coordinate (C1) at (0.5,0.8);
\coordinate (D1) at (0.5,0.5);

\fill[blue!30, opacity=0.5] (A1) -- (B1) -- (C1) -- (D1) -- cycle;

\draw[thick] (A1) -- (B1) -- (C1) -- (D1) -- cycle;

\coordinate (A2) at (0.5,0.4);
\coordinate (B2) at (0.5,0.9);
\coordinate (C2) at (0.9,0.9);
\coordinate (D2) at (0.9,0.4);

\fill[green!30, opacity=0.5] (A2) -- (B2) -- (C2) -- (D2) -- cycle;

\draw[thick] (A2) -- (B2) -- (C2) -- (D2) -- cycle;

\coordinate (A3) at (0.5,0.1);
\coordinate (B3) at (0.5,0.4);
\coordinate (C3) at (0.8,0.4);
\coordinate (D3) at (0.8,0.1);

\fill[blue!30, opacity=0.5] (A3) -- (B3) -- (C3) -- (D3) -- cycle;

\draw[thick] (A3) -- (B3) -- (C3) -- (D3) -- cycle;

\node at (0.1,0.9)    {$\Omega_0$};
\node at (0.35,0.65)  {$\Omega_1$};
\node at (0.7,0.65)  {$\Omega_2$};
\node at (0.3,0.3)    {$\Omega_3$};
\node at (0.65,0.25)  {$\Omega_4$};

\node[
  anchor=west,
  text width=4.3cm,
  align=center
] at (-0.2,1.3) {INPUT: mesh with different volumetric tags $V_1,\ldots,V_4$};

\newcommand{\Vsmall}[2]{%
\tikz[baseline=(n.base)]{
  \node[draw, fill=#2!20, rounded corners, inner sep=2pt] (n) {$V_{#1}$};
}}

\node at (0.1, -0.2)   {\Vsmall{1}{white}};
\node at (0.35, -0.2) {\Vsmall{2}{blue}};
\node at (0.6, -0.2) {\Vsmall{3}{green}};
\node at (0.85, -0.2)   {\Vsmall{4}{red}};

\node at (1.25,0.5) {$\longrightarrow$};

\begin{scope}[xshift=1.5cm]  

\draw[gray, thin] (1,0) -- (1,1);
\draw[gray, thin] (0,1) -- (1,1);

\foreach \i in {0,...,\n} {
    \foreach \j in {0,...,\n} {
        \pgfmathsetmacro\x{\i/\n}
        \pgfmathsetmacro\y{\j/\n}
        \ifnum\i<\n
            \ifnum\j<\n
                \draw[gray, thin] (\x,\y) -- ({(\i+1)/\n},\y);
                \draw[gray, thin] (\x,\y) -- (\x,{(\j+1)/\n});
                \draw[gray, thin] (\x,\y) -- ({(\i+1)/\n},{(\j+1)/\n});
            \fi
        \fi
    }
}

\coordinate (A) at (0.1,0.1);
\coordinate (B) at (0.1,0.5);
\coordinate (C) at (0.5,0.5);
\coordinate (D) at (0.5,0.1);

\draw[thick, red] (A) -- (B) -- (C) -- (D) -- cycle;

\coordinate (A1) at (0.2,0.5);
\coordinate (B1) at (0.2,0.8);
\coordinate (C1) at (0.5,0.8);
\coordinate (D1) at (0.5,0.5);

\draw[thick, green] (A1) -- (B1) -- (C1) -- (D1) -- cycle;

\coordinate (A2) at (0.5,0.4);
\coordinate (B2) at (0.5,0.9);
\coordinate (C2) at (0.9,0.9);
\coordinate (D2) at (0.9,0.4);

\draw[thick, cyan] (A2) -- (B2) -- (C2) -- (D2) -- cycle;

\coordinate (A3) at (0.5,0.1);
\coordinate (B3) at (0.5,0.4);
\coordinate (C3) at (0.8,0.4);
\coordinate (D3) at (0.8,0.1);

\draw[thick, green] (A3) -- (B3) -- (C3) -- (D3) -- cycle;

\draw[thick, black] (A1) -- (D1);
\draw[thick, black] (A3) -- (B3);

\draw[thick, purple] (C1) -- (D1);
\draw[thick, purple] (B3) -- (C3);

\draw[thick, orange] (A2) -- (0.5, 0.5);

\node[
  anchor=west,
  text width=4cm,
  align=center
] at (-0.2,1.3) {Tagging facets, so that $\Gamma_{ij}$ can be uniquely determined by tag $F$};

\node at (0.5, -0.2) {$\textcolor{green}{F_{12}}\,
\textcolor{cyan}{F_{13}}\,
\textcolor{red}{F_{14}}\,
\textcolor{purple}{F_{23}}\,
\textcolor{black}{F_{24}}\,
\textcolor{orange}{F_{34}}$};

\end{scope}

\node at (2.75,0.5) {$\longrightarrow$};

\begin{scope}[xshift=3cm]  

\node[
  anchor=west,
  text width=4.4cm,
  align=center
] at (-0.2,0.9) {Python dictionary\\ \texttt{.pickle} \vspace{0.1cm}
\\ $\Vsmall{1}{white} : \{\textcolor{green}{F_{12}},\, \textcolor{cyan}{F_{13}},\, \textcolor{red}{F_{14}}\}$
\\ $\Vsmall{2}{blue} : \{\textcolor{green}{F_{12}},\textcolor{purple}{F_{23}},\, \textcolor{black}{F_{24}}\}$
\\ $\Vsmall{3}{green} : \{\textcolor{cyan}{F_{13}}, \textcolor{purple}{F_{23}}, \textcolor{orange}{F_{34}}\}$
\\ $\Vsmall{4}{red} : \{\textcolor{red}{F_{14}},\textcolor{black}{F_{24}},\, \textcolor{orange}{F_{34}}\}$};

\draw[->] (0.5,0.4) -- ++(0,-0.2)
  node[below, align=center] {FEniCS weak forms};

\end{scope}

\end{tikzpicture}
    \end{adjustbox}
    \caption{The pre-processing starts with an input mesh where $M$ volumetric tags $V_1,\ldots,V_{M}$, with $M\leq N + 1$, are assigned to subdomains. Adjacent cells have different tags (left panel). Given this input, facets are tagged, so that each $\Gamma_{ij}$ tag is unique (middle panel). We use the following tagging strategy for facets, $F_{pq}=F_{qp}=M\cdot\min(p,q) +\max{(p,q)}$, with $p,q=1,\ldots,M$ denoting volumetric tags. Finally, we create a Python dictionary in the pickle format, mapping facet tags to volume tags (right panel), so that FEniCS weak forms can be assembled correctly.
    }
    \label{fig:chart}
\end{figure}
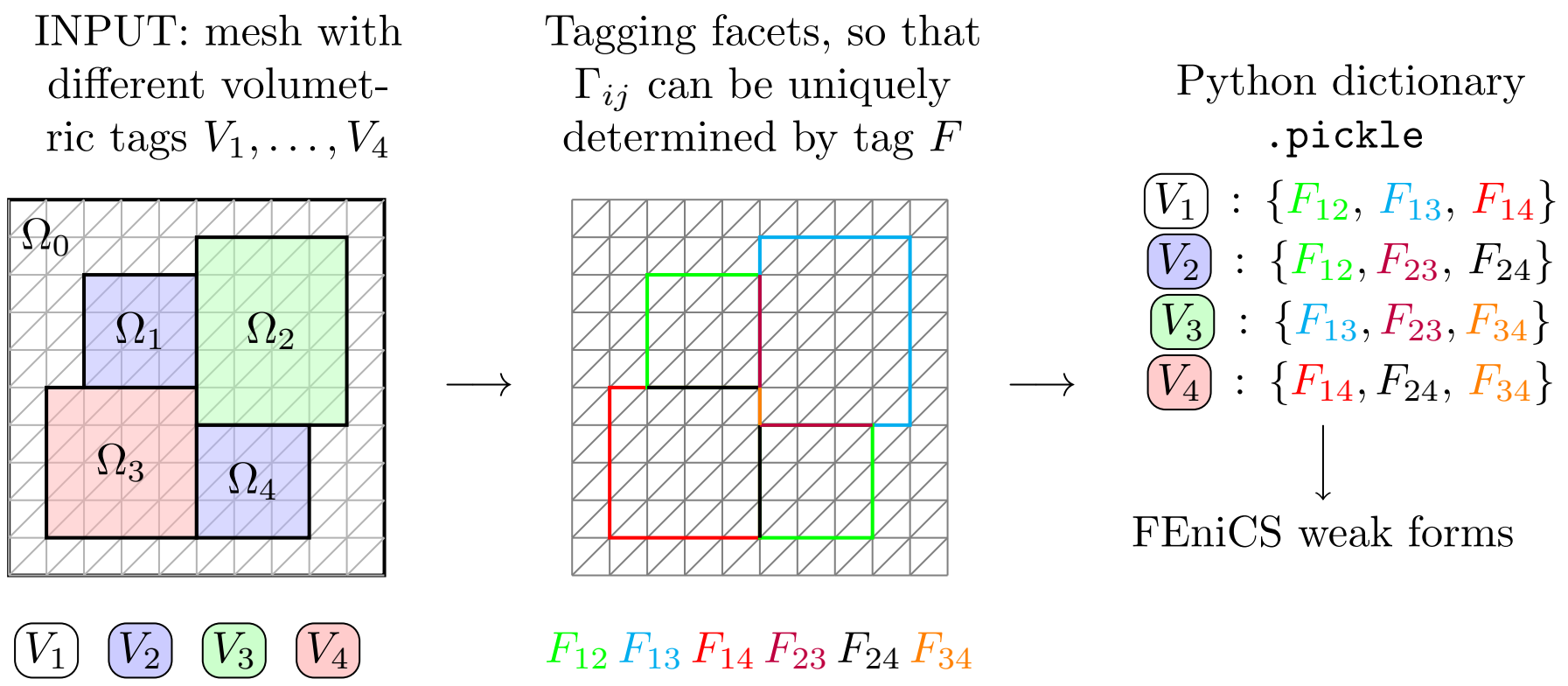

Given a mesh with different volume tags for different cells (so that cells sharing an interface have distinct tags), we construct auxiliary facet tags to mark membrane facets, along with a dictionary encoding the dependency graph between cells, which we use to guide the block assembly. In terms of implementation, disjoint cells can be grouped under the same tag, or index $i$ in Eq.~\eqref{eq:d_i}, to avoid looping over individual cells (which can be thousands). Tags, on the other hand, are heuristically bounded by $\approx5$ in realistic cardiac geometries.\footnote{No theoretical bound exists for the corresponding coloring problem in 3D.} 
From the dictionary, interface contributions for index $i$ can be isolated, so that all $\Gamma_{ij}$ terms are identified and encoded in the FEniCS weak forms following Eq.~\eqref{eq:weak_2}. 
The detailed pre-processing pipeline is illustrated and described in Figure~\ref{fig:chart}.
 
\subsection{GPU-accelerated assembly with cuDOLFINx}\label{sec:cudolfinx}
As GPU-accelerated solvers in PETSc can significantly reduce solution time compared to CPU ones~\cite{centofanti2024comparison},
the assembly of right-hand side vectors becomes a non-negligible bottleneck unless it is also accelerated. This motivates the use of cuDOLFINx to offload assembly operations to the GPU.
We remark that other (non-linear) EMI-based problems may require matrix assembly at each time step, so that assembly can become a bottleneck, see e.g.~\cite{Benedusi2024}.
As we will demonstrate in Sect.~\ref{sec:results}, CUDA assembly kernels in cuDOLFINx deliver substantial speedups over CPU implementations. 
Assembly performance is largely memory-bandwidth bound, and modern GPUs provide significantly higher memory bandwidth than CPUs. For example, on the NVIDIA GH200 used in our experiments, the GPU is equipped with HBM3 memory, providing 8x times higher memory bandwidth than the LPDDR memory available to the CPU, as confirmed by a STREAM benchmark run on our hardware. This architectural difference makes GPU acceleration particularly well-suited for assembly routines, even when accounting for multicore CPU parallelism.
Moreover, running both assembly and solve on the GPU minimizes host–device data movement, which could otherwise limit overall performance.

At assembly time, DOLFINx
iterates over all active mesh entities (cells, facets,...),
and compute the local
element tensor~$\mathcal{K}^c_e$ (or local vector~$\bm{f}^c_e$), and scatters these local contributions into the global matrix~$\mathcal{K}$ (or vector~$\bm{f}$) through a \emph{degree-of-freedom map} (dofmap).
Here and in the following, we denote the global matrix and vector of a
generic variational problem by $\mathcal{K}$ and $\bm{f}$, respectively, consistently with the notation already introduced for the EMI system in Eq.~\eqref{eq:monolithic}. 
For a bilinear form on a cell~$c$ with local
DoF indices $l,m = 0,\dots, n_{\mathrm{loc}}-1$, the scatter reads
\begin{equation}\label{eq:scatter}
  \mathcal{K}\bigl[\pi(c,l),\;\pi(c,m)\bigr]
  \;\mathrel{+}=\; \mathcal{K}^c_e[l,m],
\end{equation}
where $\pi(c,l)$ denotes the dofmap, mapping the $l$-th local
DoF of cell~$c$ to its global index.  Linear forms are handled
analogously.
This architecture separates numerical integration 
from global bookkeeping (the dofmap and the scatter).
 
Internally, cuDOLFINx reuses the C element
kernels generated by FFCx, embedding them directly within CUDA
assembly kernels.  Because GPU execution demands a much higher degree
of parallelism than CPU assembly loops, cuDOLFINx reorganizes the
data structures of DOLFINx, in particular the dofmap, mesh
coordinates, and coefficient arrays, into layouts that are
efficient for GPU memory access patterns.
After this one-time data transfer, all element integrations and
scatter operations are executed on the GPU and CPU--GPU data movement is
limited to updating coefficient and constant values between
time steps.
On a single GPU, cuDOLFINx assembly kernels can match the
throughput of FEniCS running on several many-core CPU
nodes~\cite{pachev2026cudolfinx}, and multi-GPU assembly has also
been demonstrated~\cite{pachev2025multi}.

\subsection{Restricted assembly via modified dofmaps}
For the monolithic linear system in \eqref{eq:monolithic} each contribution $\mathcal{K}_i = \tau A_i + M_i \in \mathbb{R}^{n \times n}$ has nonzero entries in the rows
corresponding to the degrees of freedom of subdomain~$\Omega_i$.
The stiffness matrix~$A_i$, arising from the volumetric integral~$a_i$, couples only DoFs within~$\Omega_i$ and thus contributes exclusively
to the diagonal block~$(i,i)$.
The mass matrix~$M_i$, arising from the interface form~$p_i$
(cf.~\eqref{eq:d_i}), contributes to both the diagonal
block, through $\int_{\Gamma_i} u_i\,\phi_i\,\mathrm{d}s$, and
to the off-diagonal blocks~$(i,j)$, through the coupling integrals
$-\int_{\Gamma_{ij}} u_j\,\phi_i\,\mathrm{d}s$ for each
neighboring subdomain~$j$.
 
To assemble this system in FEniCS, each block form is
defined on a continuous Galerkin (CG) function space that is
formally constructed on the \emph{entire} mesh but is
\emph{restricted} to the cells and facets belonging to the
relevant subdomain or interface.  The solution is therefore
continuous within each~$\Omega_i$ but discontinuous across the
interfaces~$\Gamma_{ij}$.
On CPUs, this assembly is handled by multiphenicsx,
which creates restricted function spaces and uses PETSc index sets
together with local sub-matrices, which are accessed through
\texttt{MatGetLocal\allowbreak SubMatrix} and released through
\texttt{MatRestoreLocal\allowbreak SubMatrix}, to assemble each block into the monolithic system.
For each block~$(i,j)$, multiphenicsx extracts a local sub-matrix
indexed by the block's row and column index sets, equips it with a
custom \texttt{ISLocalToGlobalMapping} that maps unrestricted DoF
indices to restricted positions (with $-1$ for inactive DoFs), and
then invokes the standard DOLFINx assembly routine on the
sub-matrix; vector assembly follows an analogous
extract--assemble--restore pattern.
Porting this strategy to GPUs would require extracting and
restoring sub-matrices across the host--device boundary for every
block, negating much of the GPUs throughput advantage.
 
Rather than replicate the PETSc sub-matrix extraction and
local-to-global remapping mechanism, cuDOLFINx
implements restricted assembly by
\emph{modifying the dofmap itself}, so that the standard GPU
assembly kernel can be reused without any changes.
The key observation is that in the scatter operation~\eqref{eq:scatter},
the only component that decides \emph{where} a local contribution
is placed in the global system is the dofmap~$\pi$.
By constructing a modified dofmap~$\tilde\pi_i$ for each subdomain
form, one can direct the contributions from form~$i$ to the
correct rows and columns of the monolithic system, while the
element integration kernel remains identical.
 
More precisely, let $V_h$ be the CG finite element space defined on
the full mesh, with $n$ global DoFs and standard dofmap~$\pi$.
For subdomain~$\Omega_i$, let
$\mathcal{R}_i \subset \{0,\dots, n - 1\}$ denote the
\emph{restriction set}: the indices of the global DoFs whose
support intersects~$\overline\Omega_i$, so that
$n_i = |\mathcal{R}_i|$ is the number of DoFs of $V(\overline\Omega_i)$
introduced in Section~\ref{sec:space_disc}.  Let
$r_i \colon \mathcal{R}_i \to \{0,\dots,n_i-1\}$ be a renumbering map
that assigns a local block index to each restricted DoF.  The monolithic offset
for subdomain~$i$ is
\begin{equation}\label{eq:offset}
  o_i = \sum_{k=0}^{i-1} n_k, \qquad o_0 = 0.
\end{equation}
The \emph{modified dofmap} for form~$i$ is then
\begin{equation}\label{eq:modified_dofmap}
  \tilde\pi_i(c,l) =
  \begin{cases}
    o_i + r_i\bigl(\pi(c,l)\bigr),
      & \text{if } \pi(c,l) \in \mathcal{R}_i, \\[3pt]
    -1, & \text{otherwise},
  \end{cases}
\end{equation}
where $-1$ is a sentinel value indicating an inactive DoF whose
contributions are discarded during the scatter.
 
With $\tilde\pi_i$ in place, the assembly of form~$i$ on the GPU
proceeds exactly as in~\eqref{eq:scatter}.
For bilinear forms of arity two (i.e., matrix assembly), the row
and column function spaces may belong to different subdomains;
accordingly, independent modified dofmaps $\tilde\pi_i^{\mathrm{row}}$
and $\tilde\pi_j^{\mathrm{col}}$ are constructed from the row and
column restriction sets, with their respective offsets.  The scatter
then reads
\begin{equation}\label{eq:restricted_scatter}
  \mathcal{K}\bigl[\tilde\pi_i^{\mathrm{row}}(c,l),\;
         \tilde\pi_j^{\mathrm{col}}(c,m)\bigr]
  \;\mathrel{+}=\; \mathcal{K}^c_e[l,m],
  \qquad c \in \mathcal{T},
\end{equation}
where $\mathcal{T}$ denotes the set of mesh entities
(cells or facets) active for the form at hand.
For vector assembly (arity one), only the test-space dofmap is
modified, and the right-hand side is assembled analogously:
\begin{equation}\label{eq:rhs_scatter}
  \bm{f}\bigl[\tilde\pi_i(c,l)\bigr]
  \;\mathrel{+}=\; \bm{f}^c_e[l],
  \qquad c \in \mathcal{T}_i.
\end{equation}
No special scatter logic, sub-matrix extraction, or host--device
copies are needed: the standard assembly kernel writes directly
into the monolithic system at the correct locations, governed
entirely by the pre-computed modified dofmap.

\subsection{Software design}
At the C++ level, the class \texttt{CUDADofMap} stores the
modified dofmap as a flat integer array on the GPU.
Its constructor receives the original DOLFINx \texttt{DofMap},
an offset~$o_i$, a separate ghost offset (for parallel ghost
DoFs, which are stored after all owned DoFs in the monolithic
vector), and an optional restriction map~$r_i$.  It ``unrolls''
the dofmap, expanding any block structure, and applies the
offset and restriction in a single pass, producing the
device-side array that the assembly kernel indexes into
(cf.~\eqref{eq:modified_dofmap}).  The class \texttt{CUDAForm}
holds a reference to the underlying DOLFINx form and, when
\texttt{set\_restriction} is called, constructs restricted
\texttt{CUDADofMap} objects for each function space of the form.
At the Python level, the \texttt{BlockCUDAForm} class orchestrates
the block assembly: it collects all subdomain forms, computes the
cumulative offsets~\eqref{eq:offset} by iterating over the
function spaces, builds the restriction sets from user-supplied
index arrays, and calls \texttt{set\_restriction} on each
\texttt{CUDAForm}.  Boundary conditions are similarly remapped
through \texttt{CUDADirichletBC}, which applies the same offsets
and restrictions to the constrained DoF indices.
 
This design offers several practical advantages.
First, the element integration kernels generated by FFCx are reused
verbatim: no kernel rewriting or recompilation is needed for
restricted assembly.
Second, the dofmap modification is a one-time setup cost; once the
modified dofmap arrays reside in GPU memory, every subsequent
assembly call (including the per-time-step right-hand side assembly)
proceeds at full GPU throughput with no additional overhead.
Third, no intermediate sub-matrices or sub-vectors are created, and
no PETSc scatter or copy operations are required; the monolithic
system is populated directly.
As a consequence, the GPU-accelerated version of the code has no multiphenicsx dependency; the restricted-assembly functionality is self-contained within cuDOLFINx.
 
\subsection{Solution strategy}\label{sec:solution_strategy}
 
To assess the performance and robustness of iterative solvers, we performed a systematic exploration of Krylov methods and preconditioning strategies. Throughout this work, we compare four solver configurations: Conjugate Gradient (CG) and Generalized Minimal Residual (GMRES), each combined with either Algebraic Multigrid (AMG) or Jacobi preconditioning.
AMG, implemented through hypre BoomerAMG, was selected as the primary preconditioning strategy because of its optimal complexity for elliptic operators. Its performance was compared against a simple Jacobi preconditioner, which is highly parallel and therefore well suited to GPU architectures. To ensure a fair comparison between different solver and preconditioner combinations, all simulations employed a relative tolerance of $10^{-7}$ and monitored convergence using the unpreconditioned residual norm.
We remark that matrix-free multigrid variants, which are likely better suited to GPUs, have recently been studied for interface problems on regular geometries~\cite{GUO2026119275}; however, no matrix-free AMG formulation currently exists, so this approach sacrifices the geometric flexibility needed for complex geometries.
 
A substantial part of this study was devoted to identifying robust AMG parameters for CPU and GPU runs. The main difference between the two architectures lies in the setup phase of AMG, namely in the algorithm that splits the unknowns into coarse (C) and fine (F) points. Unless prescribed otherwise, BoomerAMG uses HMIS coarsening with extended+i interpolation on CPUs, and PMIS coarsening with a reformulation of extended+i interpolation in terms of sparse matrix--matrix products on GPUs~\cite{Hypredocs,LiSjogreenYang2021}. Both algorithms reduce the operator complexity of classical parallel coarsening schemes~\cite{DeSterckYangHeys2006}. PMIS selects the C-points as a maximal independent set of the strength-of-connection graph, breaking ties randomly; it requires only local operations and is currently the only coarsening available on GPUs~\cite{FalgoutEtal2021,Hypredocs}. HMIS first applies the inherently sequential first pass of the classical Ruge--St\"uben algorithm within each process subdomain, and then PMIS across process boundaries~\cite{DeSterckYangHeys2006}, so that the resulting hierarchy depends on the parallel partition. Compared with HMIS, PMIS tends to produce coarser hierarchies with lower complexity, at the price of a less effective coarse-grid correction and, hence, of a slower AMG convergence~\cite{DeSterckYangHeys2006,FalgoutEtal2021,LiSjogreenYang2021}.
Depending on the geometry type (regular or realistic), we use different smoothers and coarsening types, motivated by an extensive iper-parameter search. 

For regular geometries, we used the default coarsening and interpolation of each architecture and the hybrid $L_1$-scaled SOR/Jacobi relaxation scheme\footnote{This smoother applies Gauss--Seidel within each process subdomain and a $L_1$-scaled Jacobi across subdomains. It is equivalent to an inexact block Jacobi using SOR within the blocks.}.
 In this case, the difference between the CPU and GPU preconditioners is negligible, and AMG-CG requires essentially the same number of iterations for both architectures.

For realistic geometries, whose meshes are unstructured and highly heterogeneous, the quality of the C/F splitting matters more, and the dependence of the AMG preconditioner on the parallel partition becomes relevant. Moreover, the $L_1$-scaled SOR/Jacobi smoother may yield an AMG preconditioner that is not positive definite, which can cause AMG-CG to fail. To compare the two architectures fairly, we therefore prescribe the same, partition-independent AMG configuration on CPUs and GPUs, with PMIS coarsening and extended+i interpolation, a Chebyshev polynomial smoother on all levels, including the coarsest one, and no C/F-ordered relaxation.
With these settings, CPU and GPU runs build essentially the same hierarchy 
, resulting in the same iteration counts across architectures and node counts, cf.\ Tables~\ref{tab:strong_unstructured_solver_scaling} and~\ref{tab:weak_unstructured_solver_scaling}. We chose the Chebyshev smoother over the $L_1$-scaled Jacobi one, which is also independent of the partition, because it was slightly faster per time step.
 
 In all cases, we use the default strong threshold of $0.25$, and no aggressive coarsening. The complete set of PETSc and hypre options required to reproduce the experiments is reported in Appendix~\ref{app:solver_params}.

\section{Benchmark models}\label{sec:models}
In this section, we outline the model scenarios used to evaluate the performance of the EMI framework in the cardiac setting. These models are designed to bridge the gap between idealized computational benchmarks and biologically relevant simulations, directly addressing our research questions regarding solver scalability and robustness. We employ two distinct configurations: Model A, which utilizes a modular geometry designed for scalability tests, and Model B, which represents a transition toward real-world use cases.
For clarity, we denote the number of degrees of freedom $n$ as ``DoFs'' throughout the subsequent tables.

\subsection{Model A} 
\begin{figure}
    \centering
    \includegraphics[width=0.7\linewidth]{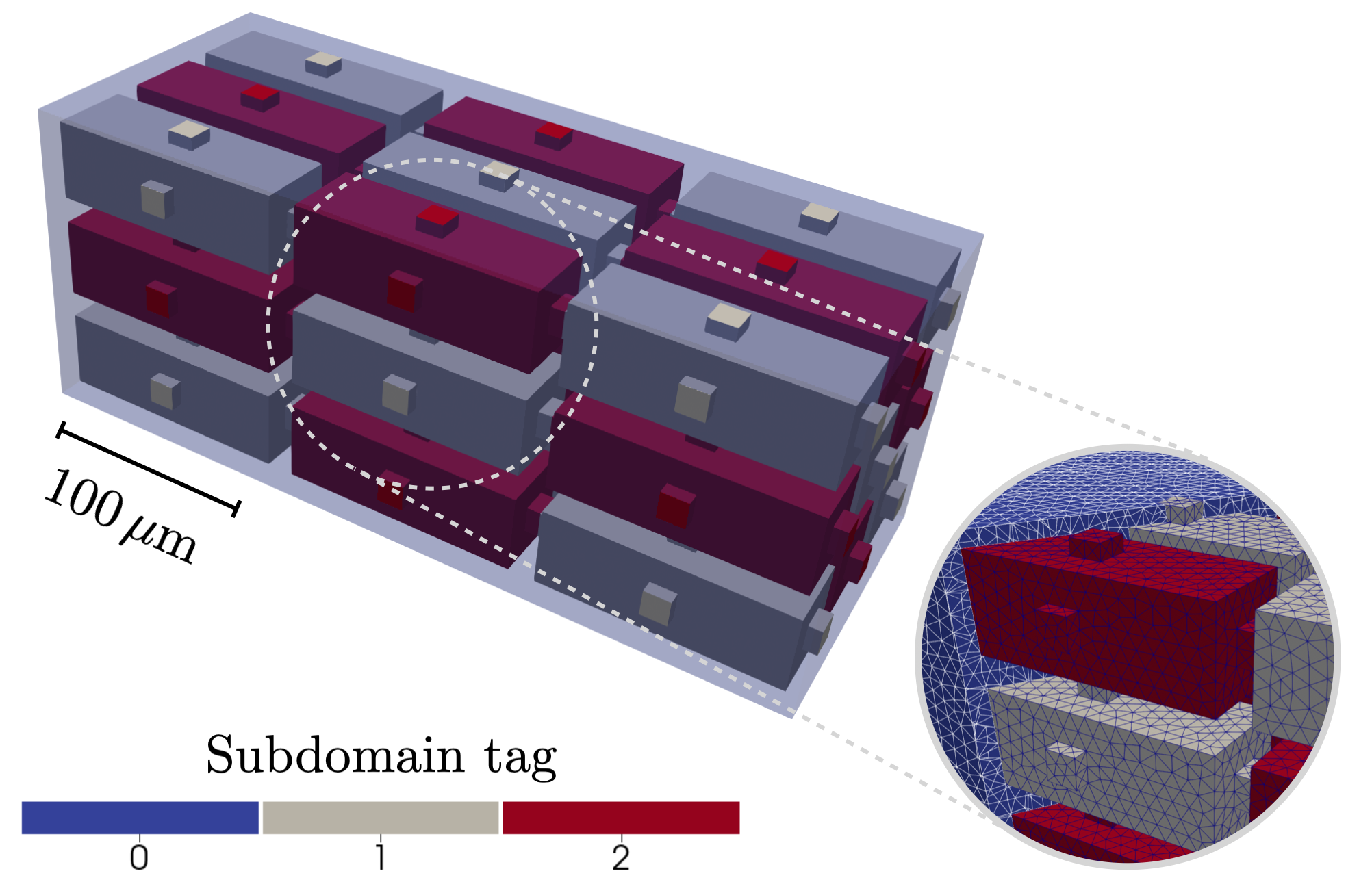}
    \caption{Model A: an example of 3$\times$3$\times$3 myocytes regular geometry. Different colors indicate the subdomain tags. This case presents a total of $M=3$ tags (see also Figure~\ref{fig:chart}).}
    \label{fig:struct}
\end{figure}

To evaluate scalability and performance in a controlled environment, we first employ a modular, synthetic geometry. In this model, cardiac cells are represented as regular blocks with smaller cubes attached to each face to represent gap junctions. This highly structured configuration is particularly suitable for performing weak and strong scalability tests; its inherent symmetry allows for an even distribution of the computational load among processors or GPUs, minimizing load imbalance and providing a clear baseline for hardware performance. 
We show an example of this geometry in Figure~\ref{fig:struct}, while the number of elements and DoFs are reported in Table~\ref{tab:mesh_struct} for each configuration.

\begin{table}[t]
    \centering
\begin{tabular}{llrr}
\toprule
Mesh name & Myocytes ($N$) & Elements & DoFs \\
\midrule
\textbf{A1} & 20$\times$10$\times$5   & 8\,290\,738 & 2\,122\,086 \\
\textbf{A2} & 20$\times$10$\times$10  & 16\,585\,884 & 4\,207\,204 \\
\textbf{A3} & 20$\times$20$\times$10  & 33\,185\,104 & 8\,377\,104 \\
\textbf{A4} & 20$\times$20$\times$20  & 66\,381\,339 & 16\,683\,928 \\
\bottomrule
\end{tabular}
\caption{Model A: number of elements and degrees of freedom (DoFs) for each configuration.}
\label{tab:mesh_struct}
\end{table}

\subsection{Model B} 
\begin{figure}
    \centering
    \includegraphics[width=0.6\linewidth]{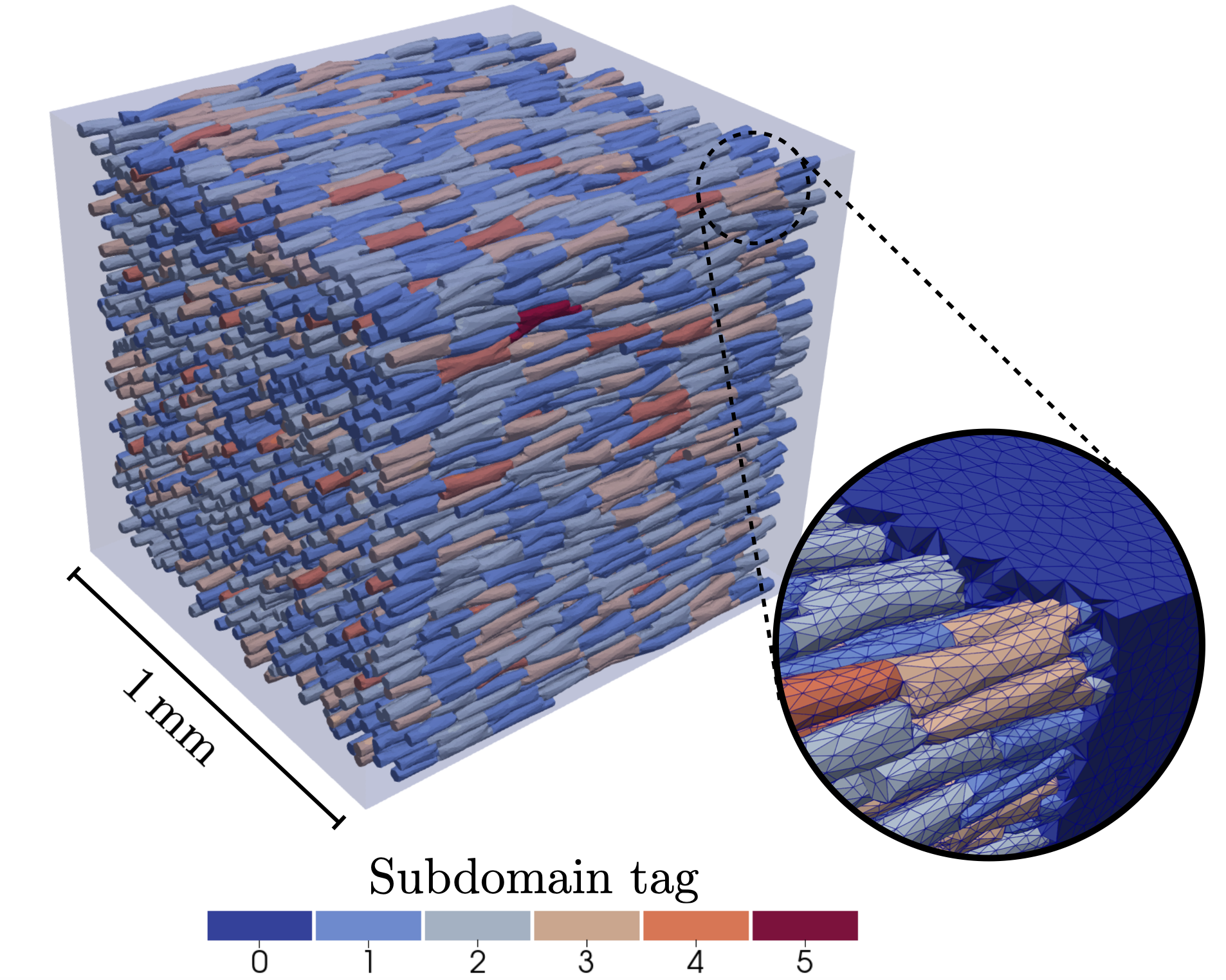}    
    \caption{Model B geometry, with the coarsest setting (B1). Different colors indicate subdomain tags. Six different tags ($M=6$) are needed to avoid neighboring myocytes with the same tag. The myocyte network is surrounded by an extracellular bath (transparent box). The magnified circle highlights the tetrahedral elements employed for the intra- and the extracellular space discretization. 
    }
    \label{fig:unstruct}
\end{figure}

\begin{figure}
    \centering
    \includegraphics[width=0.45\linewidth]{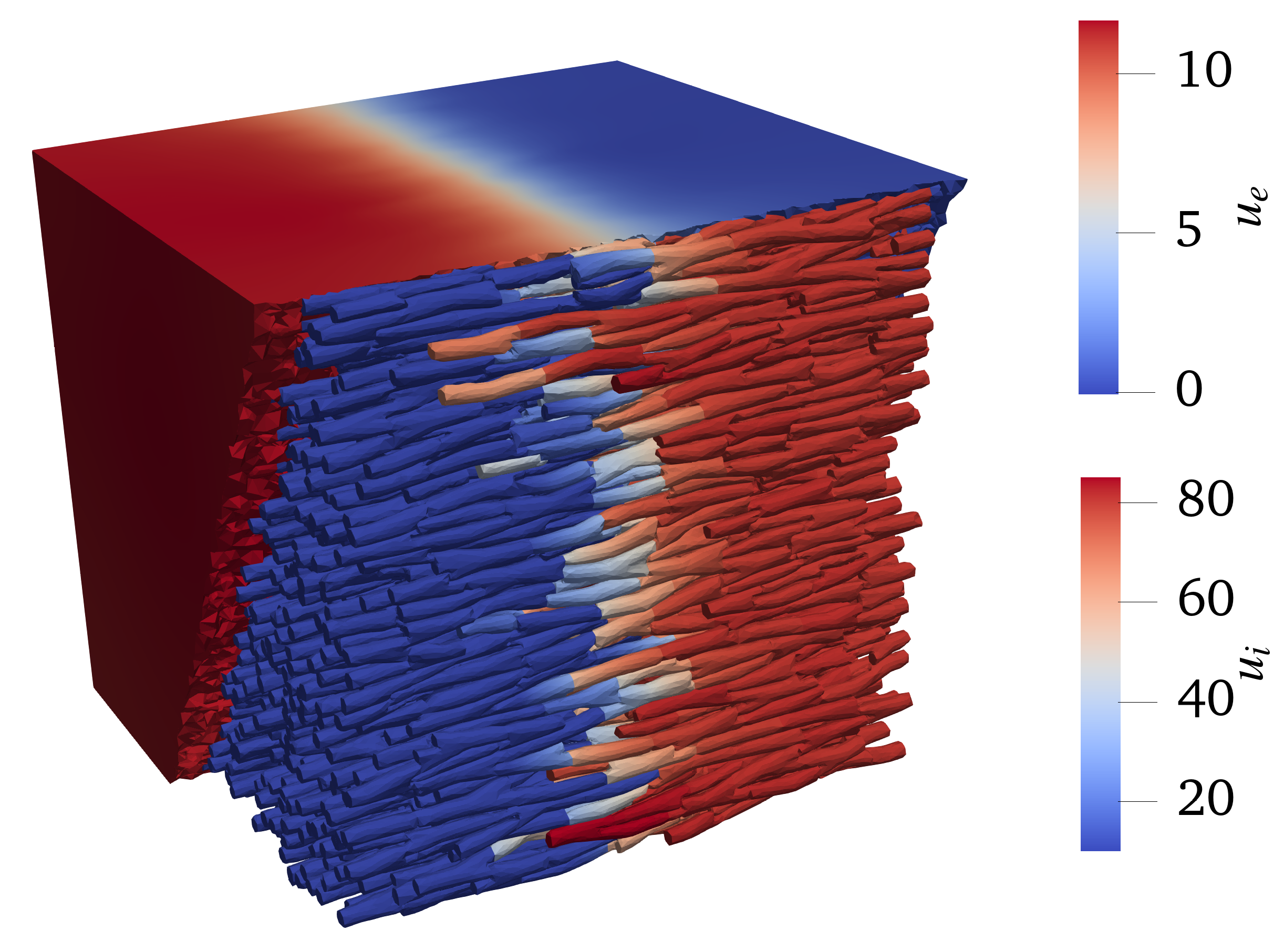}
    \includegraphics[width=0.45\linewidth]{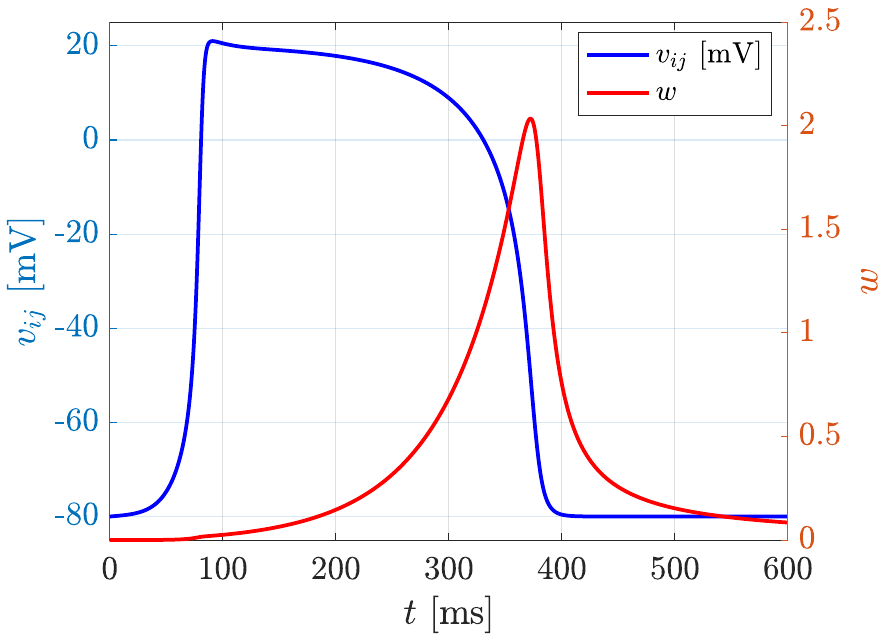}
    \caption{Model B example solution. Left panel: intracellular potential (here denoted by $u_i$) and extracellular one ($u_e$) at $t=100$ ms, for a traveling wave starting from the left part of the domain. Let us note that the potential wavefront is not homogeneous, because of the non trivial connectivity network between myocardiocytes. Right panel: corresponding membrane potential $v_{ij}=v_{01}$ and recovery variable $w$ evolution at a reference spatial point in $\Gamma_{ij}$. }
    \label{fig:sol}
\end{figure}

To assess performance in a biologically relevant scenario, we consider a second realistic geometry representing a 1\,mm$^3$
 volume of cardiac tissue. Figure~\ref{fig:unstruct} illustrates the subdomain tagging, while Figure~\ref{fig:sol} shows the corresponding solution. This geometry was created with the anni algorithm~\cite{potse-eccomas}, which is designed to synthetically generate myocyte alignments characteristic of cardiac tissue.
 In short, the algorithm first generates a point cloud of myocyte centers and links neighboring centers, splitting each link halfway. The union of the half-link domains associated with a center defines a myocyte. Heuristic constraints, such as a minimum distance between centers or a maximum myocyte length, govern the creation of this network. The network then defines a level-set function describing the cell membranes, which is discretized with the remesher \textit{mmg}~\cite{DAPOGNY2014358-mmg}.
To facilitate weak scalability tests and investigate solver robustness, the mesh is progressively refined, again with \textit{mmg}, across four levels, ensuring that surfaces are preserved without changing the topology of the myocytes alignment. The statistics of the different refinement levels are given in Table~\ref{tab:modelB_statistics}. 

We remark that this geometry is highly complex, and the meshing process can
produce nearly degenerate tetrahedra, typically where the extracellular space between neighbouring cells has almost vanishing thickness. In this work, all meshes were inspected and degenerate tetrahedra with $\rho_{\min}<1\cdot10^{-6}$ have been removed (cfr. Table~\ref{tab:modelB_statistics}). We therefore stress the importance of careful mesh generation and inspection for physiological simulations.

\begin{table}[ht]
\centering
\footnotesize
\begin{tabular}{lcccc}
\toprule
 Mesh name & \textbf{B1} & \textbf{B2} & \textbf{B3} & \textbf{B4} \\
\midrule
DoFs & $3\,391\,119$ & $6\,601\,594$ & $10\,999\,350$ & $22\,122\,183$ \\
$n_\mathrm{elems}$ & $10\,799\,476$ & $25\,447\,781$ & $44\,738\,748$ & $94\,367\,416$ \\
Myocytes ($N$) & $7\,843$ & $7\,845$ & $7\,845$ & $7\,845$ \\
$h$ & $11.7 \pm 5.6$ & $8.37 \pm 3.65$ & $6.74 \pm 2.33$ & $5.15 \pm 1.56$ \\
$\rho$ & $0.420 \pm 0.215$ & $0.436 \pm 0.263$ & $0.562 \pm 0.262$ & $0.664 \pm 0.221$ \\
$\rho_{\min}$ & $1.0 \cdot 10^{-6}$ & $1.0 \cdot 10^{-6}$ & $1.0 \cdot 10^{-6}$ & $1.0 \cdot 10^{-6}$ \\
$\Gamma_\mathrm{myo}$ & $13\,800 \pm 2\,700$ & $15\,800 \pm 3\,200$ & $16\,000 \pm 3\,200$ & $16\,100 \pm 3\,200$ \\
$V_\mathrm{ECS}/V_\mathrm{ICS}$ & $0.888$ & $0.650$ & $0.632$ & $0.617$ \\
\bottomrule
\end{tabular}
\caption{Model B statistics. If applicable, values are shown as the mean and the standard deviation over all the entities of the mesh. All spatial quantities are given in \textmu m. The variable $n_\mathrm{elems}$ denotes the number of tetrahedra in the meshes, $N$ the number of myocytes, $h$ the edge length, and $\rho$ the mesh quality, defined as the ratio of inscribed to circumscribed sphere radius (normalized to 1 for a regular tetrahedron); $\rho_{\min}$ is its minimum over the mesh. 
$\Gamma_\mathrm{myo}$ denotes the membrane area per myocyte (including cell--cell junctions) and 
$V_\mathrm{ECS}/V_\mathrm{ICS}$ the ratio between extracellular and intracellular spaces.}
\label{tab:modelB_statistics}
\end{table}

\subsection{Ionic models}
For both model A and B, we consider two separate models for the interface current $I_{ij}$, depending on the type of membrane, following equation~\eqref{eq:F_term}. Among the wide range of available ionic models, we select two models with a favorable trade-off between computational cost and physiological accuracy, capturing intercellular conduction while remaining inexpensive to evaluate. For gap junctions (i.e., interfaces between cells), we adopt a Plonsey passive linear model \cite{plonsey2007quantitative} with no additional state variables,
corresponding to purely ohmic conduction, so that only $R_g$ in \eqref{eq:F_term} is set. For the cellular membrane (i.e., interfaces involving the extracellular media), we employ the Aliev--Panfilov (AP) ionic model~\cite{aliev1996simple}, where $I_{\text{ion}}$ is a nonlinear function of the transmembrane potential and a single recovery variable $\bm{w} = w$ is present:
\begin{equation}
I_{\text{ion}}(v, w) = k v(v - a)(v - 1) + v w,
\end{equation}
where $k>0$ is an excitability parameter and $a>0$ a threshold parameter. The evolution of the gating variable $w$ is governed by
\begin{equation}
R(v, w) = \left(\epsilon + \frac{\mu_1 w}{v + \mu_2}\right)\left(-w - k v(v - a - 1)\right),
\end{equation}
where $\epsilon,\mu_1,\mu_2>0$ are recovery parameters.
The AP model is expressed in dimensionless form. In particular, following~\cite{aliev1996simple}, the normalized variable $v_{ij}$ is related to the measured transmembrane potential $E$ (in mV) through the affine transformation 
\begin{equation}
v_{ij} = \frac{E + 80}{100},
\end{equation}
so that $v_{ij} = 0$ corresponds to the resting potential $E_{\text{rest}} = -80\ \text{mV}$ and $v_{ij} = 1$ corresponds to a depolarized state with amplitude $20\ \text{mV}$. 
Similarly, time is rescaled according to
$$t_{\text{phys}} \,[\text{ms}] = \alpha_t\, t,$$
where $t$ denotes the dimensionless time variable used in the model. This scaling ensures that the action potential duration matches physiological values (the action potential duration is $ \approx 25.58$ units, while in the canine physiological case is $\approx 330\ \text{ms}$~\cite{aliev1996simple}). We provide the list of the two models parameters in Table \ref{tab:params}.

\begin{table}[t]
\centering
\begin{adjustbox}{max width=\textwidth}
\begin{tabular}{lllrl}
\toprule
 & Parameter & Symbol & Value & Unit \\
\midrule
\multirow{3}{*}{Common parameters}
& Membrane capacitance & $C_\mathrm{M}$ & 1 & \textmu F/cm$^2$ \\
& Timestep & $\Delta t$ & 0.005 & ms \\
\midrule
\multirow{1}{*}{Passive (gap junction)}
& Gap junction resistance & $R_\mathrm{g}$ & 0.003 & k$\Omega$cm$^2$ \\
\midrule
\multirow{9}{*}{Aliev--Panfilov}
& Excitability parameter & $k$ & 8.0 & -- \\
& Threshold parameter & $a$ & 0.15 & -- \\
& Recovery parameter & $\epsilon$ & 0.002 & -- \\
& Recovery dynamics & $\mu_1$ & 0.2 & -- \\
& Recovery dynamics & $\mu_2$ & 0.3 & -- \\
& Minimum potential (scaling) & $v_{ij}$ (min) & -80 & mV \\
& Maximum potential (scaling) & $v_{ij}$ (max) & 20 & mV \\
& Time scaling factor & $\alpha_t$ & 12.9 & ms$^{-1}$ \\
\bottomrule
\end{tabular}
\end{adjustbox}
\caption{Physical and model parameters. The AP model is solved in dimensionless form; $v_{ij}$ (min) and $v_{ij}$ (max) define the affine transformation between physical voltage (mV) and normalized variables, while $\alpha_t$ accounts for time rescaling between model units and milliseconds.}
\label{tab:params}
\end{table}

\section{Experimental results}
\label{sec:results}

In this section, we evaluate the performance of the GPU-accelerated implementation and compare it to the CPU version. We begin with an evaluation of the right-hand side assembly routines, and proceed to investigate the performance of the linear solver, presenting both strong and weak scaling experiments. 
The experiments for Model~A, as well as the assembly experiments for Model~B, are
run on the Vista system at the Texas Advanced Computing Center (TACC), and the
linear solver experiments for Model~B on the Alps system ``Daint'' at the Swiss
National Supercomputing Centre (CSCS). Vista provides both GPU and CPU
nodes: each GPU node hosts one NVIDIA GH200 Grace Hopper Superchip, while each
CPU node contains two 72-core NVIDIA Grace CPUs, for a total of 144 cores. On
Daint, each node is equipped with four GH200 superchips\footnote{Daint nodes are connected through an HPE Slingshot-11 network.}.
To compare CPU and GPU performance consistently across the two systems, we use
the following convention throughout the paper: a \emph{GPU node} denotes a
single GH200 driven by one MPI rank, and a \emph{CPU node} denotes 144 Grace
cores, each driven by one MPI rank.
We stress that a GPU node and a CPU node are not directly comparable in cost, power consumption, or architecture; the one-to-one node comparison used below is adopted purely for convenience in organizing the results, not as a claim of fairness.
In all the experiments, we consider one hundred time steps ($N_t=100$), and report the average results per time step. In terms of initial conditions, for all $i,j$, we consider the initial membrane potential $v^{\mathrm{in}}_{ij}=-80\cdot(1-\chi_S)$ mV, with $S\subset\Omega$ a initially activated region and $\chi_S$ its characteristic function and $w^{\mathrm{in}}_{ij}=0$.

\subsection{Right-hand side assembly}\label{sec:assembly}
We perform both strong and weak scaling tests of the right-hand side vector assembly.
We remark that matrix assembly is time-independent, so it has to be performed only once, thus it becomes negligible for large $N_t$.
The results are summarized in Table \ref{tab:strong_assembly_scaling}.  

\begin{table}[t]
    \centering
    \begin{adjustbox}{max width=\textwidth}
\begin{tabular}{lrrrrr}
\toprule
Nodes & A1 GPU & A1 CPU & B1 GPU & B1 CPU \\
\midrule
1 &  0.05 &  0.19 &  0.28 &  0.61 \\
2 &  0.03 &  0.09 &  0.15 &  0.31 \\
4 &  0.02 &  0.04 &  0.08 &  0.15 \\
8 &  0.01 &  0.02 &  0.04 &  0.08 \\
\bottomrule
\end{tabular}
\end{adjustbox}
\caption{Strong scaling of the EMI vector assembly routine for the meshes A1 and B1. Performance numbers are reported in seconds per timestep.}
\label{tab:strong_assembly_scaling}
\end{table}
The strong scaling results show that, for assembly, a single GH200 matches the performance of at least two 144-core CPU nodes.
For the weak scaling test, one node was used for A1 and B1 meshes, and the number of nodes was doubled for each refined mesh. This led to a roughly equal ratio of DoFs to nodes across runs. To facilitate comparison between runs, performance is reported in millions of DoFs per second and normalized by the number of nodes. The results are summarized in Table \ref{tab:weak_assembly_scaling}.
\begin{table}%
    \centering
    \begin{adjustbox}{max width=\textwidth}
\begin{tabular}{lrrrrr}
\toprule
Nodes & A GPU & A CPU & B GPU & B CPU\\
\midrule
1 &  40.51 &  11.02 &  11.96 &  5.54 \\
2 &  37.63 &  10.20 &  9.43 &  6.31 \\
4 &  36.45 &  9.90 &  8.87 &  6.32 \\
8 &  33.80 &  9.31 &  8.19 &  6.44 \\
\bottomrule
\end{tabular}
\end{adjustbox}
\caption{Weak scaling of the EMI vector assembly routines. Performance numbers are reported in millions of DoFs per second, normalized by the number of nodes. Larger numbers indicate better performance. Meshes A1 to A4 and B1 to B4 are used when increasing the number of nodes.} 
\label{tab:weak_assembly_scaling}
\end{table}
Unlike the strong scaling test, the weak scaling test shows different performance ratios between the GPU and CPU for the Model A and Model B meshes. For the regular meshes (A1-A4), each GPU node is the equivalent of three CPU nodes, while for the realistic meshes (B1-B4), one GPU is  approximately $30\%$ to $100\%$ more efficient than a single CPU node. 

\subsection{Linear solver}
For the regular geometry (Model A), we used AMG preconditioned CG both in the GPU and CPU cases. For GPUs, we also included results using Jacobi preconditioning, since this strategy is more effective than AMG in some cases. 
For the realistic geometry (Model B), AMG preconditioning is used in all cases, with the same configuration on CPUs and GPUs (further details in Section~\ref{sec:solution_strategy}). We perform strong and weak scaling tests, with the same mesh progression used in Section~\ref{sec:assembly} for experiment on assembly times. 

\subsubsection{Regular geometry}

For Model A, Tables~\ref{tab:strong_structured_solver_scaling} and \ref{tab:weak_structured_solver_scaling} contain strong and weak scaling results.
\begin{table}%
    \centering
    \begin{adjustbox}{max width=\textwidth}
\begin{tabular}{lrrrrr}
\toprule
Nodes & GPU Jacobi-CG & GPU AMG-CG & CPU AMG-CG \\
\midrule
1 &  0.90 (1803) &  1.41 (9) &  0.39 (10) \\
2 &  0.65 (1802) &  0.98 (10) &  0.25 (10) \\
4 &  0.52 (1802) &  0.55 (10) &  0.25 (9) \\
8 &  0.47 (1803) &  0.39 (10) &  0.44 (10) \\
\bottomrule
\end{tabular}
\end{adjustbox}
\caption{Strong scaling of the EMI linear solve routines for mesh A1. Performance numbers are reported in seconds per timestep. The number of Krylov iterations is reported in parentheses.}
\label{tab:strong_structured_solver_scaling}
\end{table}
In terms of strong scaling, we note that Jacobi preconditioning is more effective than AMG for a lower GPU count, despite the much larger number of Krylov iterations. This trend reverses as the number of GPUs nodes increases, since communication scales linearly with the number of iterations.
The CPU solver does not scale past two nodes, which can also be explained by communication overhead. Each CPU node contains 144 cores, so the CPU solver utilizes a much larger number of MPI processes than the GPU solver, which requires only one MPI process per node. Single-node GPU performance is worse than single-node CPU performance, but the GPU solver scales better.

\begin{table}%
    \centering
    \begin{adjustbox}{max width=\textwidth}
\begin{tabular}{lrrrrr}
\toprule
Nodes & GPU Jacobi-CG & GPU AMG-CG & CPU AMG-CG \\
\midrule
1 &  2.36 (1803) &  1.51 (9) &  5.38 (10) \\
2 &  1.98 (1854) &  1.79 (10) &  4.04 (10) \\
4 &  1.85 (1917) &  2.06 (10) &  2.48 (9) \\
8 &  1.51 (2265) &  2.09 (9) &  1.67 (9) \\
\bottomrule
\end{tabular}
\end{adjustbox}
\caption{Weak scaling study of linear solve routines for a sequence of regular geometry meshes (from A1 to A4 in Table~\ref{tab:mesh_struct}). Performance numbers are reported in millions of DoFs per second, normalized by the number of nodes. Larger numbers indicate better performance. Krylov iterations are reported in parentheses for the linear solvers. The runtimes in seconds are plotted in Figure~\ref{fig:structured_weak_scaling}.} 
\label{tab:weak_structured_solver_scaling}
\end{table}
The weak scaling study in  Table~\ref{tab:weak_structured_solver_scaling} confirms the same trend.  The communication overhead and the increasing iteration count explains the performance deterioration with more nodes. In terms of weak scaling, the CPU solve is sub-optimal so that the GPU implementation outperforms the CPU one on a per-node basis for the largest run (mesh A4 on 8 nodes). As expected, AMG preconditioning is robust to mesh refinement in all cases.
\begin{figure}%
    \centering
    \includegraphics[width=0.9\linewidth]{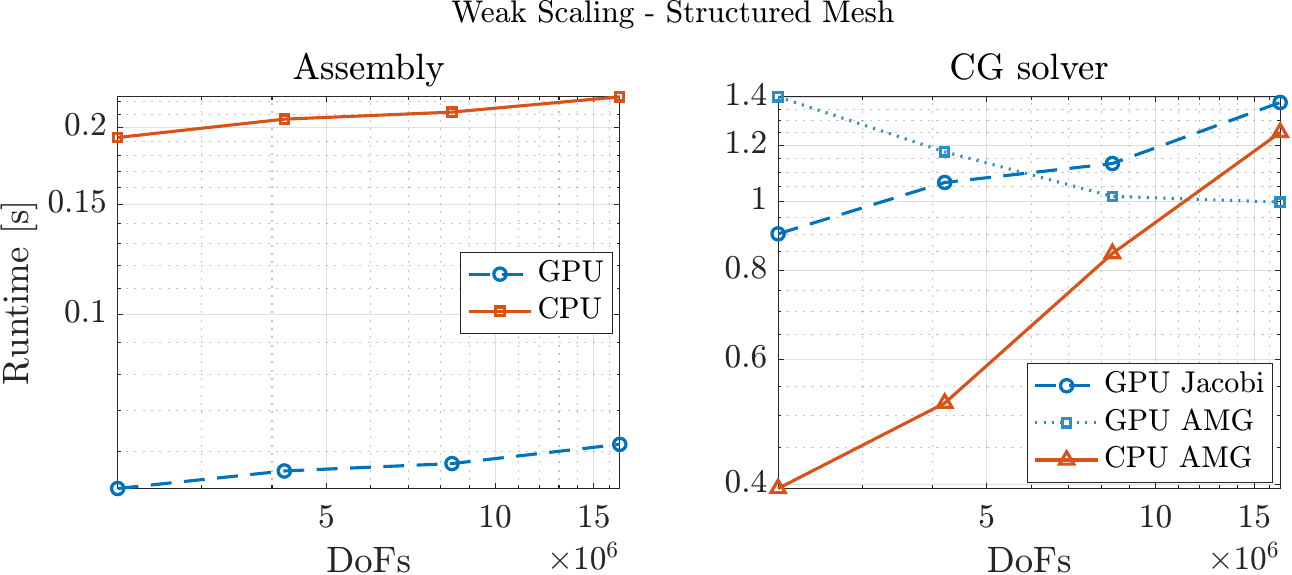}
    \caption{Average runtimes per timestep for the EMI linear solvers and assembly operations for a weak scaling test on regular meshes (A1-A4). Corresponds to Table~\ref{tab:weak_structured_solver_scaling}.
    }
    \label{fig:structured_weak_scaling}
\end{figure}
 
\subsubsection{Realistic geometry}

Tables~\ref{tab:strong_unstructured_solver_scaling} and \ref{tab:weak_unstructured_solver_scaling} contain strong and weak scaling results for Model B.
\begin{table}%
    \centering
    \begin{adjustbox}{max width=\textwidth}
\begin{tabular}{lrrrr}
\toprule
 & \multicolumn{2}{c}{B1} & \multicolumn{2}{c}{B2} \\
\cmidrule(lr){2-3} \cmidrule(lr){4-5}
Nodes & GPU & CPU & GPU & CPU \\
\midrule
1 & 0.70 (68) & 1.57 (66) & OOM       & 2.33 (57) \\
2 & 0.75 (68) & 0.99 (66) & 0.97 (56) & 1.33 (56) \\
4 & 0.67 (66) & 0.92 (66) & 0.67 (54) & 0.95 (55) \\
8 & 0.62 (68) & 0.96 (64) & 0.58 (56) & 0.87 (55) \\
\bottomrule
\end{tabular}
\end{adjustbox}
\caption{Solution strategy (AMG-CG) strong scaling for the realistic B1 and B2 meshes. Performance numbers are reported in seconds per timestep. The number of Krylov iterations is reported in parentheses. OOM denotes an out-of-memory failure.}
 \label{tab:strong_unstructured_solver_scaling}
 \end{table}
Strong scaling saturates early on B1, for both CPU and GPU. On B2 the runtimes keep decreasing up to 8 nodes, although the GPU runs out of memory on a single node. The number of iterations stays essentially constant across node counts in all cases and GPU runs are consistently faster than the corresponding CPU runs. From two to eight nodes, B2 speeds up by a factor of 1.7 on GPUs and 1.5 on CPUs.
On the GPU, the cost per Krylov iteration stays between about 9 and 12\,ms in all the B1 runs, including the single-GPU one, where no MPI communication is involved, and placing up to four GPUs on the same node made the runs slightly slower rather than faster. This floor is therefore set by the latency of the many small kernels and global reductions of each AMG V-cycle rather than by the network, and strong scaling pays off only while the work per GPU stays well above it. B2 does not fit on a single GPU: with the current implementation, the memory required by the cuDOLFINx forms limits a GH200 (96\,GB) to about five million DoFs.
\begin{table}%
    \centering
    \begin{adjustbox}{max width=\textwidth}
\begin{tabular}{lrrrrr}
\toprule
Model & GPU AMG-CG & GPU AMG-GMRES & CPU AMG-CG & CPU AMG-GMRES \\
\midrule
B1 & 4.81 (68) & 4.32 (72) & 2.17 (66) & 1.89 (71) \\
B2 & 3.40 (56) & 3.23 (56) & 2.48 (56) & 2.57 (51) \\
B3 & 3.80 (48) & 2.91 (60) & 2.54 (48) & 2.01 (58) \\
B4 & 1.94 (50) & 1.76 (52) & 1.75 (58) & 1.42 (69) \\
\bottomrule
\end{tabular}
\end{adjustbox}
\caption{Linear solver weak scaling for a sequence of realistic geometry meshes. Performance numbers are reported in millions of DoFs per second, normalized by the number of nodes. Larger numbers indicate better performance. Krylov iteration counts are reported in parentheses for the linear solvers. The corresponding runtimes in seconds are depicted in Figure~\ref{fig:unstructured_weak_scaling}.}
\label{tab:weak_unstructured_solver_scaling}
\end{table}
Table~\ref{tab:weak_unstructured_solver_scaling} reports the weak scaling of the linear solvers. The number of Krylov iterations generally decreases as the mesh is refined. This is consistent with the theory in~\cite{benedusi2024dense}: uniform refinement shifts the problem, spectrally, closer to a pure Laplacian, for which robust or even faster convergence is expected. The iteration counts therefore do not explain the loss of per-node throughput observed at larger sizes. On the GPU the AMG-CG throughput drops from 4.81 to 1.94 MDoFs/s per node between B1 and B4, while on the CPU it stays roughly constant up to B3 and then falls to 1.75. This degradation is likely due to parallel overheads at higher node counts. The GPU is faster than the CPU on B1--B3, up to about a factor two on B1, while the two are comparable on B4. AMG-CG is generally as fast as or faster than AMG-GMRES, with very similar iteration counts.
We remark that, with the configuration described in Section~\ref{sec:solution_strategy} and the meshes cleaned of nearly degenerate tetrahedra, a single set of solver options converged for all the entries of Tables~\ref{tab:strong_unstructured_solver_scaling} and~\ref{tab:weak_unstructured_solver_scaling}, with no per-entry tuning.

\begin{figure}[H]
    \centering
    \includegraphics[width=0.9\linewidth]{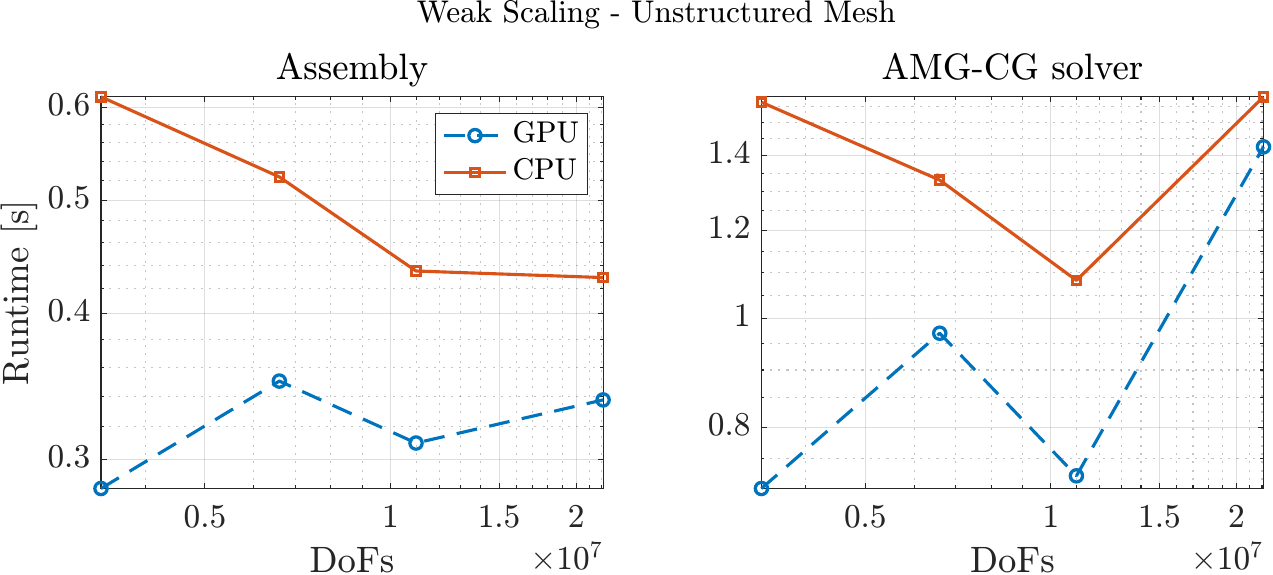}
    \caption{Average runtimes per timestep for the EMI AMG-CG linear solver and assembly operations for a weak scaling test on realistic geometry meshes. Corresponds to Table~\ref{tab:weak_unstructured_solver_scaling}.
    }
    \label{fig:unstructured_weak_scaling}
\end{figure}

 \section{Conclusions}
\label{sec:conclusion}

This work presents the first performance evaluation and comparison of GPU-accelerated linear algebra routines within cuDOLFINx. We consider the EMI (Extracellular-Membrane-Intracellular) model as a mixed-dimensional benchmark on unstructured meshes. EMI is an elliptic model posed on a domain partitioned into subdomains coupled through interface conditions, and its numerical solution is of significant practical interest ~\cite{tveito2021modeling,benedusi2024modeling,centofanti2025parallel,SteyerBiosig}. The presence of interfaces makes the EMI problem more challenging to assemble and solve than a standard elliptic problem, making it a meaningful stress test for the cuDOLFINx framework.

We presented a GPU-accelerated implementation of the EMI model and evaluated its performance against a CPU-based implementation across both regular and realistic geometries. We assessed right-hand side assembly and linear solver performance using strong and weak scaling experiments on multiple super-computing systems. 
For the assembly routines, strong scaling results showed that a single GH200 GPU node is roughly equivalent to two 144-core CPU nodes. 
Weak scaling results revealed a similar but geometry-dependent trend: each GPU node was equivalent to approximately three CPU nodes for the regular meshes, and $30\%$ to $100\%$ more efficient than a single CPU node for the realistic meshes. In all cases, the weak scaling of assembly routines was close-to-optimal.
For linear solvers, the results were more nuanced and depended on geometry, preconditioning strategy, and problem scale. For the regular geometry, Jacobi preconditioning outperformed AMG on GPUs at low node counts despite requiring substantially more Krylov iterations, though this trend reversed as the number of nodes increased and communication costs began to dominate. The CPU solver failed to scale past two nodes, a consequence of the much larger number of MPI processes per node compared to the single-process-per-node GPU configuration. In the weak scaling regime, the GPU implementation exhibits near-optimal scaling and  eventually overtook the CPU solver on a per-node basis at the largest problem size tested.
For realistic meshes, both architectures are robust in terms of iteration counts and show similar strong scaling, with the runtimes close to saturation after a few nodes. A GPU node is consistently faster than a CPU node, by a factor between 1.3 and 2.2. In the weak scaling tests, the CPU solver retains a close-to-optimal parallel efficiency, whereas the GPU solver scales less well but remains overall faster, with one GPU node being between 1.1 and 2.2 times faster than one CPU node. The iteration counts are robust also in this case.
We emphasize that the one-to-one node comparisons used throughout this study are adopted for convenience and do not reflect equivalence in cost, power consumption, or hardware architecture between GPU and CPU nodes. 
Overall, the experiments show, for both assembly and solution, that one GH200 GPU node is equivalent to roughly 160 to 310 CPU cores, depending on the setting.

Finally, we remark that robust convergence on the realistic geometries requires care. Hybrid smoothers, which apply Gauss--Seidel within each process subdomain and Jacobi across subdomains, depend on the parallel partition and may yield AMG preconditioners that are not positive definite. Removing degenerate mesh elements after careful inspection, and adopting the same partition-independent AMG configuration, based on Chebyshev smoothing, on both architectures, made a single solver setting converge in all the tests, with iteration counts essentially independent of the architecture and of the number of nodes. On the GPU, the strong scalability on these meshes is limited by the latency of the AMG cycle and by the memory footprint of the current implementation (about five million DoFs per GH200), which are natural targets for future work.

\appendix

\section{Detailed solver parameters}
\label{app:solver_params}

This appendix lists the specific configuration parameters used for the PETSc solvers and hypre BoomerAMG preconditioners. For a detailed explanation of these options, the reader is referred to~\cite{Hypredocs, Petscdocs}. We used PETSc 3.23.7, and hypre 2.32.0.

\subsection{Regular geometry configurations (Model A)}
The solver configuration for regular geometries employed the following baseline AMG settings both for CPU and GPU:

\begin{tcolorbox}
\begin{verbatim}
# -ksp_type here is cg or gmres
-ksp_rtol: 1e-7
-ksp_norm_type: unpreconditioned
-pc_hypre_boomeramg_relax_type_all: l1scaled-SOR/Jacobi
-pc_hypre_boomeramg_agg_nl: 0
-pc_hypre_boomeramg_strong_threshold: 0.25
\end{verbatim}
\end{tcolorbox}
Coarsening and interpolation were left to the defaults of each architecture (HMIS on CPUs, PMIS on GPUs, see Section~\ref{sec:solution_strategy}).
For comparison, the Jacobi preconditioner was configured as follows:

\begin{tcolorbox}
\begin{verbatim}
-ksp_type: cg
-ksp_rtol: 1e-7
-ksp_norm_type: unpreconditioned
-pc_type: jacobi
\end{verbatim}
\end{tcolorbox}

\subsection{Realistic geometry configurations (Model B)}
For the realistic models, CPU and GPU runs used exactly the same options. CG, with homogeneous Neumann conditions and the constant nullspace, and GMRES, with one Dirichlet point, share the following configuration:
\begin{tcolorbox}
\begin{verbatim}
# -ksp_type here is cg or gmres
-ksp_rtol: 1e-7
-ksp_norm_type: unpreconditioned
-pc_type: hypre
-pc_hypre_boomeramg_coarsen_type: PMIS
-pc_hypre_boomeramg_interp_type: ext+i
-pc_hypre_boomeramg_relax_type_all: Chebyshev
-pc_hypre_boomeramg_relax_type_coarse: Chebyshev
-pc_hypre_boomeramg_no_CF: true
-pc_hypre_boomeramg_agg_nl: 0
-pc_hypre_boomeramg_strong_threshold: 0.25
-pc_hypre_boomeramg_min_coarse_size: 2
\end{verbatim}
\end{tcolorbox}
PMIS coarsening and extended+i interpolation are the GPU defaults; they are set explicitly so that the CPU runs use them as well. C/F-ordered relaxation, enabled by default on CPUs only, is disabled. Aggressive coarsening is disabled because PETSc fixes different interpolation operators on the aggressively coarsened levels on CPUs (multipass) and GPUs (two-stage extended+e). The minimum coarse-grid size of two prevents a $1\times1$ coarsest matrix, which vanishes for the singular Neumann operator. The hypre statistics (number of levels and operator complexity) confirmed that CPU and GPU runs build the same hierarchy. 

\section*{Acknowledgments}
BP gratefully acknowledges the use of the Vista system at the Texas Advanced Computing Center under the ``ADCIRC" allocation. 
EC acknowledges the CINECA award under the ISCRA initiative, for the availability of high-performance computing resources and support (project DDO2CARD).
EC, JS and JT have been supported by the European High-Performance Computing Joint Undertaking (JU) EuroHPC under grants agreement No 955495 (MICROCARD) and No 101172576 (MICROCARD-2). The JU receives support from the European Union’s Horizon 2020 research and innovation program and France, Italy, Germany, Austria, Norway, Switzerland.
PB and EC acknowledge the support by grants from the Swiss National Supercomputing Centre (CSCS) under project IDs sm111 and lp157 on Alps. 
The authors would like to thank Joshua Steyer’s PI, Axel Loewe, for helpful discussions and for the freedom to pursue this independent study. The authors also wish to thank Mark Potse for his help with the mesh utilities. AI tools were used for language and grammar editing; the authors assume responsibility for all content.

\newpage
\bibliographystyle{siamplain}
\bibliography{references}

\begin{thebibliography}{10}

\bibitem{aliev1996simple}
{\sc R.~R. Aliev and A.~V. Panfilov}, {\em A simple two-variable model of
  cardiac excitation}, Chaos, Solitons \& Fractals, 7 (1996), pp.~293--301.

\bibitem{multiphenicsx}
{\sc F.~Ballarin}, {\em {multiphenicsx} - easy prototyping of multiphysics
  problems in {FEniCSx}}, 2020, \url{https://multiphenics.github.io/}.

\bibitem{BarattaEtal2023}
{\sc I.~A. Baratta, J.~P. Dean, J.~S. Dokken, M.~Habera, J.~S. Hale, C.~N.
  Richardson, M.~E. Rognes, M.~W. Scroggs, N.~Sime, and G.~N. Wells}, {\em
  {DOLFINx: The next generation FEniCS problem solving environment }}, Dec.
  2023, \url{https://doi.org/10.5281/zenodo.10447666}.

\bibitem{Benedusi2024}
{\sc P.~Benedusi, A.~J. Ellingsrud, H.~Herlyng, and M.~E. Rognes}, {\em
  Scalable approximation and solvers for ionic electrodiffusion in cellular
  geometries}, SIAM Journal on Scientific Computing, 46 (2024), pp.~B725--B751,
  \url{https://doi.org/10.1137/24M1644717}.

\bibitem{benedusi2024dense}
{\sc P.~Benedusi, P.~Ferrari, M.~Causemann, and S.~Serra-Capizzano}, {\em Dense
  cell-by-cell systems of {PDEs}: approximation, spectral analysis, and
  preconditioning}, BIT Numerical Mathematics, 66 (2026),
  \url{https://doi.org/10.1007/s10543-026-01136-6}.

\bibitem{benedusi2024modeling}
{\sc P.~Benedusi, P.~Ferrari, M.~E. Rognes, and S.~Serra-Capizzano}, {\em
  Modeling excitable cells with the {EMI} equations: spectral analysis and
  iterative solution strategy}, Journal of Scientific Computing, 98 (2024),
  p.~58, \url{https://doi.org/10.1007/s10915-023-02449-2}.

\bibitem{causemann2026pulsating}
{\sc M.~Causemann, V.~Vinje, and M.~E. Rognes}, {\em The pulsating brain: An
  interface-coupled fluid-poroelastic interaction model of the cranial},
  Mathematical Modelling of the Human Brain II: From Glymphatics to Deep
  Learning,  (2026), p.~47.

\bibitem{centofanti2025parallel}
{\sc E.~Centofanti, N.~M.~M. Huynh, L.~F. Pavarino, and S.~Scacchi}, {\em
  Parallel algebraic multigrid solvers for composite discontinuous galerkin
  discretization of the cardiac {EMI} model in heterogeneous media}, Computer
  Methods in Applied Mechanics and Engineering, 442 (2025), p.~118001.

\bibitem{centofantigdsw}
{\sc E.~Centofanti, N.~M.~M. Huynh, L.~F. Pavarino, and S.~Scacchi}, {\em {GDSW
  Preconditioners for the EMI Model in Cardiac Electrophysiology}}, in
  Proceedings of Domain Decomposition Methods in Science and Engineering XXIX,
  Springer, 2026, pp.~249--256.

\bibitem{centofanti2024comparison}
{\sc E.~Centofanti and S.~Scacchi}, {\em {A comparison of Algebraic Multigrid
  Bidomain solvers on hybrid {CPU-GPU} architectures}}, Computer Methods in
  Applied Mechanics and Engineering, 423 (2024), p.~116875.

\bibitem{DAPOGNY2014358-mmg}
{\sc C.~Dapogny, C.~Dobrzynski, and P.~Frey}, {\em Three-dimensional adaptive
  domain remeshing, implicit domain meshing, and applications to free and
  moving boundary problems}, Journal of Computational Physics, 262 (2014),
  pp.~358--378,
  \url{https://doi.org/https://doi.org/10.1016/j.jcp.2014.01.005}.

\bibitem{DeSterckYangHeys2006}
{\sc H.~De~Sterck, U.~M. Yang, and J.~J. Heys}, {\em Reducing complexity in
  parallel algebraic multigrid preconditioners}, SIAM Journal on Matrix
  Analysis and Applications, 27 (2006), pp.~1019--1039,
  \url{https://doi.org/10.1137/040615729}.

\bibitem{Ellingsrud2025}
{\sc A.~J. Ellingsrud, P.~Benedusi, and M.~Kuchta}, {\em A splitting,
  discontinuous {G}alerkin solver for the cell-by-cell electroneutral
  {Nernst–Planck} framework}, SIAM Journal on Scientific Computing, 47
  (2025), pp.~B477--B504, \url{https://doi.org/10.1137/24M1653367}.

\bibitem{FalgoutEtal2021}
{\sc R.~D. Falgout, R.~Li, B.~Sj{\"o}green, L.~Wang, and U.~M. Yang}, {\em
  Porting \textit{hypre} to heterogeneous computer architectures: {S}trategies
  and experiences}, Parallel Computing, 108 (2021), p.~102840,
  \url{https://doi.org/10.1016/j.parco.2021.102840}.

\bibitem{favino2025fractures}
{\sc M.~Favino}, {\em {Fractures and thin heterogeneities as Robin-Wentzell
  interface conditions}}, Applied Mathematical Modelling,  (2025), p.~116355.

\bibitem{fevola2021automated}
{\sc E.~Fevola, F.~Ballarin, L.~Jimenez-Juan, P.~Triverio, S.~Grivet-Talocia,
  and G.~Rozza}, {\em Automated estimation of patient-specific boundary
  conditions for cardiovascular simulations: an optimal control approach}, in
  2021 Summer Biomechanics, Bioengineering, \& Biotransport Conference, SB3C
  Foundation, Inc., 2021.

\bibitem{franzone2014mathematical}
{\sc P.~C. Franzone, L.~F. Pavarino, and S.~Scacchi}, {\em Mathematical cardiac
  electrophysiology}, Springer, 2014.

\bibitem{fumagalli2024discontinuous}
{\sc I.~Fumagalli}, {\em Discontinuous {G}alerkin method for a
  three-dimensional coupled fluid-poroelastic model with applications to brain
  fluid mechanics}, Mathematics in Engineering, 7 (2025), pp.~130--161.

\bibitem{gobel2026bddc}
{\sc F.~G{\"o}bel, N.~M.~M. Huynh, F.~Chegini, L.~F. Pavarino, M.~Weiser,
  S.~Scacchi, and H.~Anzt}, {\em A {BDDC} preconditioner for the cardiac {EMI}
  model in three dimensions}, SIAM Journal on Scientific Computing, 48 (2026),
  pp.~A646--A667.

\bibitem{GUO2026119275}
{\sc Z.~Guo, C.~Cui, and Q.~Zhang}, {\em {GPU}-based matrix-free multigrid
  methods of stable generalized finite element methods for interface problems},
  Computer Methods in Applied Mechanics and Engineering, 461 (2026), p.~119275,
  \url{https://doi.org/https://doi.org/10.1016/j.cma.2026.119275}.

\bibitem{heinlein2018frosch}
{\sc A.~Heinlein, A.~Klawonn, S.~Rajamanickam, and O.~Rheinbach}, {\em {FROSch:
  A fast and robust overlapping Schwarz domain decomposition preconditioner
  based on Xpetra in Trilinos}}, in International Conference on Domain
  Decomposition Methods, Springer, 2018, pp.~176--184.

\bibitem{huynh2023convergence}
{\sc N.~M.~M. Huynh, F.~Chegini, L.~F. Pavarino, M.~Weiser, and S.~Scacchi},
  {\em Convergence analysis of {BDDC} preconditioners for composite {DG}
  discretizations of the cardiac cell-by-cell model}, SIAM Journal on
  Scientific Computing, 45 (2023), pp.~A2836--A2857.

\bibitem{Hypredocs}
{\em Hypre online documentation}.
\newblock \url{https://hypre.readthedocs.io/en/latest/solvers-boomeramg.html}.
\newblock Accessed: 2026-03-27.

\bibitem{jaeger2021efficient}
{\sc K.~H. J{\ae}ger, K.~G. Hustad, X.~Cai, and A.~Tveito}, {\em Efficient
  numerical solution of the {EMI model representing the extracellular space
  (E), cell membrane (M) and intracellular space (I) of a collection of cardiac
  cells}}, Frontiers in Physics, 8 (2021), p.~579461.

\bibitem{LiSjogreenYang2021}
{\sc R.~Li, B.~Sj{\"o}green, and U.~M. Yang}, {\em A new class of {AMG}
  interpolation methods based on matrix-matrix multiplications}, SIAM Journal
  on Scientific Computing, 43 (2021), pp.~S540--S564,
  \url{https://doi.org/10.1137/20M134931X}.

\bibitem{pachev2025multi}
{\sc B.~Pachev, J.~Trotter, and I.~Baratta}, {\em {Multi-GPU acceleration of
  FEniCSx codes with cuDOLFINxx}}, in FEniCS 2025, 2025.

\bibitem{pachev2026cudolfinx}
{\sc B.~A. Pachev, J.~D. Trotter, and I.~A. Baratta}, {\em {cuDOLFINx}: A
  {CUDA} extension for {FEniCSx}}, in The FEniCS Project: The FEniCS 2024
  Conference, Springer, 2026, pp.~27--38.

\bibitem{Petscdocs}
{\em {PETSc} online documentation}.
\newblock \url{https://petsc.org/release/manualpages/}.
\newblock Accessed: 2026-03-27.

\bibitem{plonsey2007quantitative}
{\sc R.~Plonsey, R.~C. Barr, and A.~Bioelectricity}, {\em Quantitative
  approach}, Cham, Switzerland: Springer,  (2007).

\bibitem{potse-eccomas}
{\sc M.~Potse, L.~Cirrottola, and A.~Froehly}, {\em A practical algorithm to
  build geometric models of cardiac muscle tissue}, in 8th European Congress on
  Computational Methods in Applied Sciences and Engineering (ECCOMAS 2022), 01
  2022, \url{https://doi.org/10.23967/eccomas.2022.027}.

\bibitem{shapson2024petavoxel}
{\sc A.~Shapson-Coe, M.~Januszewski, D.~R. Berger, A.~Pope, Y.~Wu, T.~Blakely,
  R.~L. Schalek, P.~H. Li, S.~Wang, J.~Maitin-Shepard, et~al.}, {\em A
  petavoxel fragment of human cerebral cortex reconstructed at nanoscale
  resolution}, Science, 384 (2024), p.~eadk4858.

\bibitem{SteyerCheginiPotseetal.2023}
{\sc J.~Steyer, F.~Chegini, M.~Potse, A.~Loewe, and M.~Weiser}, {\em Continuity
  of microscopic cardiac conduction in a computational cell-by-cell model}, in
  2023 Computing in Cardiology Conference (CinC), vol.~50, 2023,
  \url{https://doi.org/10.22489/CinC.2023.385}.

\bibitem{SteyerBiosig}
{\sc J.~Steyer, F.~Chegini, T.~Stary, M.~Potse, M.~Weiser, and A.~Loewe}, {\em
  Electrograms in a cardiac cell-by-cell model}, in Workshop Biosignals,
  G{\"o}ttingen, Germany, Feb. 2024,
  \url{https://doi.org/10.47952/gro-publ-194},
  \url{https://doi.org/10.47952/gro-publ-194}.

\bibitem{tveito2021modeling}
{\sc A.~Tveito, K.-A. Mardal, and M.~E. Rognes}, {\em Modeling excitable
  tissue: The EMI framework}, Springer, 2021.

\end{thebibliography}

\end{document}